\documentclass[journal,twoside,web]{ieeecolor}
\usepackage{generic}
\usepackage{cite}
\usepackage{amsmath,amssymb,amsfonts}
\usepackage{algorithmic}
\usepackage{graphicx}
\usepackage{algorithm,algorithmic}
\usepackage{hyperref}
\hypersetup{hidelinks=true}  
\usepackage{color}
\usepackage{verbatim}
\usepackage{amsmath}
\usepackage{amssymb}
\makeatletter
\usepackage{cuted}
\usepackage{bbm}
\usepackage{multicol}

\renewcommand{\t}{^{\mbox{\tiny {T}}}}
\newcommand{\eproof}{\hfill\rule{2mm}{2mm}}
\newcommand{\bstate}{\begin{state} }
\newcommand{\estate}{ \hfill  \rule{1mm}{2mm}\end{state}}
\newcommand{\bass}{\begin{ass} }
\newcommand{\eass}{\hfill\rule{1mm}{2mm}\end{ass}}
\newcommand{\bpro}{\begin{property}}
\newcommand{\epro}{\hfill\rule{1mm}{2mm}\end{property}}
\newcommand{\brem}{ \begin{remark}  }
\newcommand{\erem}{\hfill \rule{1mm}{2mm}
\end{remark} }
\newcommand{\bthm}{\begin{theorem}  }
\newcommand{\ethm}{ \hfill  \rule{1mm}{2mm}
\end{theorem} }
\newcommand{\blem}{\begin{lemma}  }
\newcommand{\elem}{ \hfill \rule{1mm}{2mm}
\end{lemma} }
\newcommand{\bcorollary}{\begin{corollary}  }
\newcommand{\ecorollary}{  \hfill \rule{1mm}{2mm}
\end{corollary} }
\newcommand{\bdefn}{\begin{definition}}
\newcommand{\edefn}{  \hfill \rule{1mm}{2mm}
\end{definition} }
\newcommand{\bproposition}{\begin{proposition} }
\newcommand{\eproposition}{\hfill \rule{1mm}{2mm}
\end{proposition} }
\newcommand{\bexample}{\begin{example} \rm}
\newcommand{\eexample}{ \hfill \rule{1mm}{2mm}
\end{example} }
\newcommand{\proofnow}{\noindent{\bf Proof: }}
\newtheorem{theorem}{\bf Theorem}[section]
\newtheorem{ass}{\bf Assumption}[section]
\newtheorem{lemma}{\bf Lemma}[section]
\newtheorem{definition}{\bf Definition}[section]
\newtheorem{remark}{\bf Remark}[section]
\newtheorem{corollary}{\bf Corollary}[section]
\newtheorem{proposition}{\bf Proposition}[section]
\newtheorem{example}{\bf Example}[section]

\newtheorem{property}{\bf Property}[section]
\allowdisplaybreaks[4]
\newcommand{\prooflater}[1]{\noindent{\bf Proof of #1: }}
\newcommand{\sint}{\textstyle{\int}}
\newcommand{\ssum}{\textstyle{\sum}}
\makeatother
\usepackage{babel}

\begin{document}
\title{Achieving Distributed  Convex Optimization within Prescribed Time For High-Order Nonlinear Multi-Agent Systems }

\author{Gewei Zuo, Lijun Zhu, Yujuan Wang, Zhiyong Chen and Yongduan Song, \IEEEmembership{Fellow, IEEE}
\thanks{This work was supported in part by the National Natural Science Foundation of China under Grant 62173155 and Grant 52188102; in part by the Program for Huazhong University of Science and Technology (HUST) Academic Frontier Youth Team; and in part by the Taihu Lake Innovation Fund for Future Technology, HUST. \emph{(Corresponding Author: Lijun Zhu})}
\thanks{Gewei Zuo and Lijun Zhu are with School of Artificial Intelligence and Automation, Huazhong University of Science and Technology, Wuhan 430072, China. Lijun Zhu is also with Key Laboratory of Imaging Processing and Intelligence Control, Huazhong University of Science and Technology, Wuhan 430074, China (Emails: gwzuo@hust.edu.cn;  ljzhu@hust.edu.cn).}
\thanks { Yujuan Wang and Yongduan Song are with  School of Automation, Chongqing University, Chongqing, 400044, China (Email: yjwang66@cqu.edu.cn; ydsong@cqu.edu.cn).}
\thanks {Zhiyong Chen is with School of Engineering, The University of Newcastle, Callaghan, NSW 2308, Australia (Email: chen@newcastle. edu.cn). }
}

\maketitle
\begin{abstract}
 This paper addresses the distributed prescribed-time convex optimization (DPTCO) problem for high-order nonlinear multi-agent systems (MASs) under undirected connected graphs. A cascade design framework is proposed that divides the DPTCO implementation into distributed optimal trajectory generator design and local reference trajectory tracking controller design. The DPTCO problem is then transformed into the prescribed-time stabilization problem of a cascaded system. Using changing Lyapunov functions and time-varying state transformations with sufficient conditions, we establish criteria for prescribed-time stabilization and prove the boundedness of internal signals in closed-loop MASs. The framework addresses robust DPTCO for chain-integrator MASs with disturbances through the introduction of novel sliding-mode variables and time-varying gains. It also solves adaptive DPTCO for strict-feedback MASs with parameter uncertainty via backstepping method and descending power state transformation. Two numerical examples verify the theoretical results.
\end{abstract} 


\section{Introduction}
Distributed convex optimization (DCO) has gained widespread attention due to its numerous applications in multi-agent systems (MASs). These applications include reliable communications in wireless networks, collision avoidance among multiple robots, economic dispatch in power systems, distributed optimal power flow, traffic management for large-scale railway networks, and traffic metering in urban street networks.
In a typical DCO scenario, each agent has a local objective function that is private to itself. The global objective function is the sum of these local functions. The goal is to design distributed controllers that, using only limited local information, enable the output or state of each agent to converge to the optimum of the global objective function.  The earliest work on DCO can be tracked back to \cite{john38}, and it has attracted increasing interest in the last decade after the pioneer works in \cite{nedic2009distributed6}.

Research in DCO focuses on four key areas: broadening the scope of objective functions \cite{kia2015distributed18,rahili2016distributed17,yuan2015regularized50,le2017neurodynamic51} and systems \cite{gharesifard2013distributed23,li2022exponential11,zhang2017distributed4,huang2019distributed2,shi2022finite52,tang2020optimal3}, enhancing convergence rates \cite{ning2017distributed27,liu2023multi14,chen2018fixed8,shi2022finite52,wang2020distributed20,lin2016distributed12}, and improving disturbance rejection \cite{pilloni2020sliding5,zhang2017distributed4,wang2020distributed21,tang2020optimal3,wang2017distributed55}.
Optimization algorithms have been devised for time-invariant objective functions \cite{kia2015distributed18}, time-varying ones \cite{rahili2016distributed17}, and those with constraints \cite{yuan2015regularized50,le2017neurodynamic51}.  In \cite{li2017distributed15},
the convexity of local objective function and strong convexity of global
objective function are respectively removed. Some works aim to achieve
DCO for more general systems, such as single-integrator systems in
\cite{gharesifard2013distributed23}, linear systems in \cite{li2022exponential11},
Euler-Lagrange systems in \cite{zhang2017distributed4}, and strict-feedback
systems in \cite{huang2019distributed2,qin2022adaptive16}. Using sliding-mode
control and backstepping methods, the DCO controller can handle systems
that are high-order and nonlinear \cite{tang2020optimal3}. A common
approach to solving  DCO for high-order systems is the cascade
design where the solution to the DCO problem is divided into two parts.
The first is distributed optimum seeking, which by utilizing 
local information interaction generates local optimal references for
each agent that asymptotically converge to the optimum of
the global objective function. The second  is to design local tracking
controllers to make the output or state asymptotically converge to
the local optimal references.

Convergence rate and disturbance rejection are pivotal concerns in DCO. Studies such as \cite{lin2016distributed12,shi2022finite52} explore finite-time convergence, where agents reach consensus within a finite time frame while optimizing the global objective function. The finite-time DCO for chain integrator MASs facing mismatched disturbances is addressed in \cite{wang2020distributed21}. Furthermore, fixed-time convergence, wherein the finite settling time does not depend on initial conditions, is demonstrated in \cite{ning2017distributed27,wang2020distributed20,chen2018fixed8}. In \cite{liu2023multi14}, predefined-time DCO is accomplished by designing specific time-based functions, allowing solutions to converge to a neighborhood of the optimum within a specified period and to the exact optimum as time approaches infinity. However, this approach is not extendable to scenarios involving disturbances and high-order systems.

In this paper, we address the DPTCO for high-order nonlinear MASs with uncertainties, for which
the solution converges to the optimum within any prescribed time.
The prescribed-time control is proposed to ensure that the settling
time does not depend on the initial values and control parameters
\cite{song2017time36,zuo2022adaptive35,2024zhang57}. The main contributions of
this paper are summarized as follows.

First, a DFTCO framework for a class of nonlinear MASs with disturbances
is proposed. By embedding a cascade design, the DFTCO implementation
is divided into two parts, namely, distributed optimal trajectory
generator design and local reference trajectory tracking controller
design. The DPTCO problem is then transformed into the prescribed-time
stabilization problem of two cascaded subsystems. The first one
is for the error of the distributed estimation towards the global
optimum, and the second one is for local tracking errors. Changing
Lyapunov function and time-varying state methods, along with
some sufficient conditions, are proposed to prove the prescribed-time
stabilization of the cascaded system as well as the uniform boundedness
of internal signals in the closed-loop system. A specific distributed
optimal trajectory generator is constructed to show that the distributed
estimation errors converge to zero within a prescribed time.

Second, under the DPTCO framework, we propose a robust DPTCO algorithm
for a class of nonlinear chain-integrator MASs with external disturbances.
We design a novel sliding-mode variable and introduce a new time-varying state
transformation, which converts the prescribed-time stabilization problem
of local tracking error and other states unrelated to the output into
the boundedness of the new variable. Unlike traditional sliding-mode
control in \cite{edwards1998sliding39} and the prescribed-time work
in \cite{song2017time36}, our approach does not require the high-order
derivatives of the reference trajectory for tracking. Moreover, our
proposed algorithm is robust against any bounded external disturbances.

Third, we consider adaptive DPTCO for a class of strict-feedback MASs
with parameter uncertainty. We introduce a time-varying state transformation
of a descending power to compensate for the growth of increasing rate
induced by the derivative of time-varying gains in recursive steps. The
backstepping method with prescribed-time dynamic filters is adopted
to avoid the utilization of high-order derivatives of reference trajectory,
and an adaptive law is designed to compensate for parameter uncertainty.

The rest of the paper is organized as follows. Section \ref{sec:nominal_solution}
presents the notation and problem formulation. Section \ref{sec:A-Cascade-Design}
introduces the DPTCO framework for a type of nonlinear systems, for
which Section \ref{sec:Optimal-Trajectory-Generator} elaborates on the
optimal trajectory generator design. Given the DPTCO framework and
optimal trajectory generator, robust DPTCO for chain-integrator MASs
and adaptive DPTCO for strict-feedback MASs are discussed in Sections
\ref{sec:DPTCO-for-Chain} and \ref{sec:Adaptive-DPTCO-for}, respectively.
Numerical simulation is conducted in Section \ref{sec:Simulation-Results}, 
and the paper is concluded in Section \ref{sec:Conclusion}.
\section{Notations and Problem Formulation \label{sec:nominal_solution}}
\subsection{Notations}
$\mathbb{R}$, $\mathbb{R}_{\geq0}$, and $\mathbb{R}^{n}$ denote
the set of real numbers, the set of non-negative real numbers, and
the $n$-dimensional Euclidean space, respectively. $t_{0}$ denotes
the initial time, $T$ the prescribed-time scale, and $\mathcal{T}_{p}:=\{t:t_{0}\leq t<T+t_{0}\}$
the corresponding time interval.  The symbol $1_{N}\in\mathbb{R}^{N}$ (or $0_{N}\in\mathbb{R}^{N}$)
denotes an $N$-dimensional column vector whose elements are all $1$
(or $0$). For $\alpha\in\mathcal{K}_{\infty}$, $(\alpha(s))^{-1}=1/\alpha(s)$
for $s\in\mathbb{R}_{\geq0}/0$, while $\alpha^{-1}(s)$ is the inverse
function of $\alpha(s)$ for $s\in\mathbb{R}_{\geq0}$.

An undirected graph is denoted as $\mathcal{G}=(\mathcal{V},\mathcal{E})$,
where $\mathcal{V}=\{1,\cdots,N\}$ is the node set and $\mathcal{E}\subseteq\mathcal{V}\times\mathcal{V}$
is the edge set. The existence of an edge $(i,j)\in\mathcal{E}$ means
that nodes $i$ and $j$ can communicate with each other. Denote by $\mathcal{A}=[a_{ij}]\in\mathbb{R}^{N\times N}$
the weighted adjacency matrix, where $(j,i)\in\mathcal{E}\Leftrightarrow a_{ij}>0$
and $a_{ij}=0$ otherwise. A self-edge is not allowed, i.e., $a_{ii}=0$.
The Laplacian matrix $\mathcal{L}$ of graph $\mathcal{G}$ is denoted
as $\mathcal{L}=[l_{ij}]\in\mathbb{R}^{N\times N}$, where $l_{ii}=\sum_{j=1}^{N}a_{ij}$, and 
$l_{ij}=-a_{ij}$ for $i\neq j$. If $\mathcal{G}$ is connected,
then the null space of $\mathcal{L}$ is spanned by $1_{N}$, and
all the other $N-1$ eigenvalues of $\mathcal{L}$ are strictly positive.

\subsection{Problem Formulation}
Consider the nonlinear MASs
\begin{equation}
\begin{aligned}\dot{x}^{i} & =g_{x}^{i}(x^{i},y^{i},u^{i},d^{i}(t)),\\
\dot{y}^{i} & =g_{y}^{i}(x^{i},y^{i},u^{i},d^{i}(t)),\;i\in\mathcal{V},
\end{aligned}
\label{eq:nominal_sys}
\end{equation}
where $x^{i}\in\mathbb{R}^{n}$, $y^{i}\in\mathbb{R}^{m}$, $u^{i}\in\mathbb{R}^{q}$
are the system state, output, and control input of the$i$-th agent, respectively.
$d^{i}:[t_{0},\infty)\to\mathbb{D}\subset\mathbb{R}^{n_{d}}$
denotes the system's uncertainties or external disturbances, where
$\mathbb{D}$ is a compact set belonging to $\mathbb{R}^{n_{d}}$
and is possibly time-varying. $g_{x}^{i}:\mathbb{R}^{n}\times\mathbb{R}^{m}\times\mathbb{R}^{q}\times\mathbb{D}\to\mathbb{R}^{n}$,
$g_{y}^{i}:\mathbb{R}^{n}\times\mathbb{R}^{m}\times\mathbb{R}^{q}\times\mathbb{D}\to\mathbb{R}^{m}$
are smooth functions of their arguments, satisfying $g_{x}^{i}(0,y^{i},0,d^{i})=0$
and $g_{y}^{i}(0,y^{i},0,d^{i})=0$ for any $y^{i}\in\mathbb{R}^{m}$
and $d^{i}\in\mathbb{R}^{n_{d}}$. The output feedback
system (\ref{eq:nominal_sys}) encompasses various specific types \cite{chen2005global41},
such as  chain-integrator systems \cite{song2017time36}, strict-feedback
systems \cite{krstic1992adaptive42}, and feedforward systems \cite{xudong2003universal43}. 

In this paper, we consider the following convex optimization problem 
\begin{equation}
 \min_{y}\;\ssum_{i=1}^{N}f^{i}(y^{i}),\quad 
  \mbox{s.t. }\; \;y^{i}=y^{j},\; \forall i\neq j,
\label{eq:opti_problem}
\end{equation}
where $y=[y_{1}\t,\cdots,y_{N}\t]\t$ is the lumped output of MASs in \eqref{eq:nominal_sys}, and $f^i(y^i), i\in \mathcal V$ is the local scalar objective function, which is convex and known only to agent $i$. Motivated by the results in \cite{zhang2017distributed4}, this paper assumes that gradient function $\nabla f^i(\cdot)$ of local objective function is available. Due to equality constraints, the optimum $y^*$ of optimization problem \eqref{eq:opti_problem} has the form $y^* = 1_N \otimes z^*$ for some $z^*\in \mathbb R^m$.  

The objective of the DPTCO
is, for any prescribed time $T>0$, using local information interactions to design distributed controllers
$u^{i},i\in \mathcal V$ such that the  outputs $y$ converge to the optimum
$y^{*}$ within $T+t_0$, 
i.e.,
\begin{equation}
\lim_{t\to T+t_{0}} y(t)-y^{*}=0\label{eq:objective}
\end{equation}
irrespective of system's initial value and any other control parameters
besides $T$. Moreover, the state $x^{i}$, the output $y^{i}$, and
control input $u^{i}$ must be bounded, i.e.,
$
\Vert[(x^{i}(t))\t,(y^{i}(t))\t,(u^{i}(t))\t]\Vert<\infty
$
holds for $i\in\mathcal{V}$ and $t\in\mathcal{T}_{p}$.

A few discussions are needed to distinguish the distributed finite-time,
fixed-time, predefined-time, and prescribed-time convex optimization.
Let us call $t_{0}+T$ in the objective (\ref{eq:objective}) the
settling time. When the objective (\ref{eq:objective}) is modified
to $y^{i}(t)=z^{*},\forall t\geq T(x(t_{0}),c),i\in\mathcal{V}$ where
the function $T(x(t_{0}),c)$ depends on the initial condition $x(t_{0})$
and the parameter vector $c=[c^{1};\cdots,c^{N}]\in\mathbb{R}^{Nm_{c}}$
with $c^{i}\in\mathbb{R}^{m_{c}}$ being the control parameters for
$i$th agent, the problem is called distributed finite-time convex
optimization. In other words, the settling time $t_{0}+T$ is upper
bounded by $T(x(t_{0}),c)$ \cite{lin2016distributed12,shi2022finite52,wang2020distributed21}.
When the upper bound of the settling time is rendered to be independent
of the initial condition, i.e.,  $T+t_{0}\leq T(c)$, then the problem
is called distributed fixed-time convex optimization \cite{ning2017distributed27,wang2020distributed20,chen2018fixed8}.
If the settling time satisfies $T+t_{0}\leq T^{*}$ with a given $T^{*}$,
then the problem is called distributed predefined-time convex optimization
\cite{jimenez2020lyapunov37}. In comparison to predefined-time control
where the upper bound of settling time can be prescribed, $T+t_{0}$
is the exact settling time of system in prescribed-time control \cite{song2017time36,zuo2022adaptive35,2024zhang57}. In
the predefined-time control for MASs, the actual settling time $T(x^{i}(t_{0}),c^{i},T^{*})$
of agent $i$ is affected by initial values and control parameters.
Therefore, even though $T^{*}$ is given a priori for the whole MASs,
the actual settling time is inconsistent since initial values and
control parameters of each subsystem might be different.

In order to achieve the DPTCO,
the function
\begin{equation}
\mu(t)=1/(T+t_{0}-t)\label{eq:mu}
\end{equation}
is used throughout the paper as the time-varying gain. The function
$\mu(t)$ increases to infinity as $t$ approaches the prescribed time
$T+t_{0}$ and is commonly used in prescribed-time control. For
$t\in\mathcal{T}_{p}$, one has $\mu\in\mathbb{R}_{p}:=[1/T,\infty).$
We simplify $\mu(t)$ as $\mu$ throughout this paper if no confusion
occurs. For any $\alpha\in\mathcal{K}_{\infty}$ and $\iota\in\mathbb{R}$,
define
\begin{equation}
\begin{gathered}
\kappa^{\iota}\left(\alpha(\mu)\right)=\exp\big(\iota\sint_{t_{0}}^{t}\alpha(\mu(\tau))\mathrm{d}\tau\big),\quad t\geq t_0,\\
\kappa^{\iota}\left(\alpha(\mu(\tau))\right)=\exp\big(\iota\sint_{t_{0}}^{\tau}\alpha(\mu(s))\mathrm{d}s\big),\quad \tau\geq t_0,
\end{gathered}\label{eq:kappa_t}
\end{equation}
where we note that  $\kappa^{\iota}\left(\alpha(\mu)\right)$ converges to
zero as $t\to T+t_{0}$ for any $\iota<0$ and $\alpha\in\mathcal{K}_{\infty}$.
We study the problem under these two common assumptions.

\bass\label{ass:graph} The undirected graph $\mathcal{G}$ is connected.\eass

\bass\label{ass:cost_func} For each $i\mathcal{\in V},$
the function $f^{i}$ is first-order differentiable, and $f^{i}$
as well as its gradient $\nabla f^{i}$ are only
known to $i$-th agent. Moreover, it is $\rho_{c}$-strongly convex
and has $\varrho_{c}$-Lipschitz gradients, i.e., for $x,y\in\mathbb{R}^{m}$, $(\nabla f^{i}(x)-\nabla f^{i}(y))\t(x-y)  \geq\rho_{c}\Vert x-y\Vert^{2}$ and $\Vert\nabla f^{i}(x)-\nabla f^{i}(y)\Vert  \leq\varrho_{c}\Vert x-y\Vert$, 
where $\rho_{c}$ and $\varrho_{c}$ are positive constants.\eass

Under Assumption \ref{ass:cost_func}, $f$ is strongly 
convex as $f^{i}$ is for $i\in\mathcal{V}$. Therefore, if the optimization
 problem (\ref{eq:opti_problem})
is solvable, the optimum is unique. We need the following assumption
for the optimization problem to be sensible.

\bass \label{ass:solvable} The optimal value of global objection function (\ref{eq:opti_problem}),
denoted as $f^{*}$, is finite and the optimum set 
\begin{equation}
Y_{\mbox{\small{opt} }}=\{Y= 1_N \otimes z^*\mid\ssum_{i=1}^{N}f^{i}(z^*)=f^{*}\} \label{eq:y_opt}
\end{equation}
is nonempty and compact. \eass

\bdefn \label{def:KL} 
A function $\alpha:[0,\infty)\to[a,\infty)$
is said to belong to class $\mathcal{K}_{\infty}^{e}$, it is strictly
increasing and $\alpha(0)=a\geq0$.
A continuous function $\beta:[0,c)\times[0,T)\to[0,\infty)$
is said to belong to class $\mathcal{KL}_{T}^{e}$ if, for each fixed
$s$, the mapping $\beta(r,s)$ belongs to class $\mathcal{K}_{\infty}^{e}$
with respect to $r$ and, for each fixed $r$, the mapping $\beta(r,s)$
is decreasing with respect to $s$ and satisfies $\beta(r,s)\to0$
as $s\to T$. The function $\beta$ is said to belong to class $\mathcal{KL}_{T}^ {}$
if $\beta$ belongs to class $\mathcal{KL}_{T}^{e}$ and for each
fixed $s$, the mapping $\beta(r,s)$ belongs to class $\mathcal{K}_{\infty}$
with respect to $r$. \edefn

\bdefn \cite{holloway2019prescribed45} Consider
the system $\dot{\chi}=g(t,\chi,d(t))$ where $\chi\in\mathbb{R}^{n}$
is the state and $d(t):[t_{0},\infty)\to\mathbb{R}^{n_{d}}$ is
the external input. For any given $T>0$, the $\chi$-dynamics is
said to be prescribed-time stable if there exists $\beta\in\mathcal{KL}_{T}^{e}$
such that for $\chi_{0}\in\mathbb{R}^{n}$ and $d\in\mathbb{R}^{n_{d}}$,
$\|\chi(t)\|\leq\beta(\|\chi_{0}\|,t-t_{0})$ holds for $t\in\mathcal{T}_{p}$
where $\chi_{0}=\chi(t_{0})$. \edefn

\bdefn The continuously differentiable function $V(x):\mathbb{R}^{n}\to\mathbb{R}_{\geq0}$
is called the prescribed-time stable Lyapunov function for the system
$\dot{x}=f(x,\mu)$, if $V(x)$ and its derivative along the trajectory
of the system satisfy, for all $x\in\mathbb{R}^{n}$ and $t\in\mathcal{T}_{p}$,
\begin{equation}
\underline{\alpha}(\Vert x\Vert)\leq V(x)\leq\bar{\alpha}(\Vert x\Vert),\quad \dot{V}\leq-\tilde{\alpha}(\mu)V,\label{eq:PTLF}
\end{equation}
where $\underline{\alpha}$, $\bar{\alpha}$, $\tilde{\alpha}$ are
$\mathcal{K}_{\infty}$ functions and $\mu$ is defined in (\ref{eq:mu}). The term
$\tilde{\alpha}(\mu)$ is referred to as the prescribed-time convergent gain.
The inequalities in (\ref{eq:PTLF}) are simplified as $V(x)\sim\{\underline{\alpha},\bar{\alpha},\tilde{\alpha}|\dot{x}=f(x,\mu)\}$.
The continuously differentiable function $V(x):\mathbb{R}^{n}\to\mathbb{R}_{\geq0}$
is called the prescribed-time input-to-state stable (ISS) Lyapunov
function for the system $\dot{x}=f(x,d,\mu)$ with $d\in\mathbb{R}^{n_{d}}$
being the external input, if $V(x)$ and its derivative along the
trajectory of the system satisfy, for all $x\in\mathbb{R}^{n}$ and
$t\in\mathcal{T}_{p}$,
\begin{equation}
\begin{gathered}\underline{\alpha}(\Vert x\Vert)\leq V(x)\leq\bar{\alpha}(\Vert x\Vert),\\
\dot{V}\leq-\tilde{\alpha}(\mu)V+\tilde{\sigma}(\mu)\sigma(\Vert d\Vert)
\end{gathered}
\label{eq:ISS-PTLF}
\end{equation}
with $\underline{\alpha}$, $\bar{\alpha}$, $\tilde{\alpha}$, $\tilde{\sigma}\in\mathcal{K}_{\infty}$
and $\sigma\in\mathcal{K}_{\infty}^{e}$. $\tilde{\alpha}(\mu)$,
$\tilde{\sigma}(\mu)$ and $\sigma(\Vert d\Vert)$ are called prescribed-time
convergent, prescribed-time ISS gain and (normal)
ISS gain, respectively. The inequalities in (\ref{eq:ISS-PTLF}) are
simplified as $V(x)\sim\{\underline{\alpha},\bar{\alpha},\tilde{\alpha},[\sigma,\tilde{\sigma}]|\dot{x}=f(x,d,\mu)\}$.
When $d$ contains multiple inputs as $d=[d_{1}\t,\cdots,d_{n}\t]\t$
where $d_{i}\in\mathbb{R}^{n_{i}}$, the second inequality of (\ref{eq:ISS-PTLF})
becomes $\dot{V}\leq-\tilde{\alpha}(\mu)V+\sum_{i=1}^{n}\tilde{\sigma}_{i}(\mu)\sigma_{i}(\|d_{i}\|),$
and the inequalities are simplified as $V(x)\sim\{\underline{\alpha},\bar{\alpha},\tilde{\alpha},[\sigma_{1},\tilde{\sigma}_{1}],\cdots,[\sigma_{n},\tilde{\sigma}_{n}]|\dot{x}=f(x,d,\mu)\}$.
\edefn

\bdefn
For functions $ x_{1}(s) : \mathbb{R}_{\geq 0} \to \mathbb{R}_{\geq 0} $ and $ x_{2}(c, s) : \mathbb{R}_{\geq 0} \times \mathbb{R}_{\geq 0} \to \mathbb{R}_{\geq 0} $, we denote $ x_{1}(s) = \mathcal{S}[x_{2}(c, s)] $ to mean that, for any  $ c <\infty$, the following holds 
$\sup_{s \in \mathbb{R}_{\geq 0}} \left[ {x_{1}(s)}/{x_{2}(c, s)} \right] < \infty$. 
\edefn

\section{A Cascade Design Approach \label{sec:A-Cascade-Design}}
The cascade design approach has been used for the distributed convex
optimization problem in \cite{huang2019distributed2,tang2020optimal3,zhang2017distributed4}.
Following the cascade design principle, the optimal agreement can
be decomposed into two subproblems, namely distributed optimum
seeking and local reference trajectory tracking. To this end, we propose
the controller in the general form of
\begin{gather}
\dot{\zeta}^{i}=h_{\zeta}^{i}\left(\zeta^{i},\chi^{i},\mu\right),\quad \dot{\varsigma}^{i}=h_{\varsigma}^{i}(\varsigma^{i},\zeta^{i},\xi^{i},\mu)\label{eq:zeta}\\
u^{i}=h_{u}^{i}\left(\varsigma^{i},\zeta^{i},\xi^{i},\mu\right),\;i\in\mathcal{V},\label{eq:con}
\end{gather}
where $\chi^{i}=\Sigma_{j\in\mathcal{N}_{i}}a_{ij}(\zeta^{j}-\zeta^{i})$
is the relative information received by $i$-th agent from its neighbors
and $\xi^{i}=[(x^{i})\t,(y^{i})\t]\t$. $\zeta^{i}$-dynamics is designed
to estimate $z^{*}$ in \eqref{eq:y_opt}.
The state of $\zeta^{i}$ can be decomposed as $\zeta^{i}=[(\varpi^{i})\t,(p^{i})\t]\t$
where $p^{i}$-dynamics can be designed to adaptively find the gradient
of the local objective function $\nabla f^{i}(z^{*})$. $\zeta^{i}$-dynamics
is similar to a PI controller and designed to admit the equilibrium
point $\zeta^{i}=\zeta^{*,i}:=[(z^{*})\t,(p_{0}^{i}(z^{*}))\t]\t$
with some known function $p_{0}^{i}$. $\varsigma^{i}\in\mathbb{R}^{m_{\varsigma}}$
is the local controller state used to construct the actual control
input $u^{i}$ for tracking.

\subsection{Coordinate Transformation and Cascaded Error System}

For $i\in\mathcal{V}$, define the error states
\begin{equation}
\begin{gathered}e_{\varpi}^{i}=\varpi^{i}-z^{*},\;e_{p}^{i}=p^{i}-p_{0}^{i}(z^{*}),\;e^{i}=y-z^{*},\\
e_{y}^{i}=y^{i}-\varpi^{i},\;e_{s}^{i}=\left[(x^{i})\t,(e_{y}^{i})\t,(\varsigma^{i})\t\right]\t.
\end{gathered}
\label{eq:e_y_varpi}
\end{equation}
Note that $e_{\varpi}^{i}$ and $e_{p}^{i}$ are the errors from the
distributed optimal value seeking, $e^{i}$ is the optimal tracking
error and $e_{y}^{i}$ is the local tracking error towards the local
estimated optimal value $\varpi^{i}$. Define the lumped vectors $e_{\varpi}=\left[(e_{\varpi}^{1})\t,\cdots,(e_{\varpi}^{N})\t\right]\t$,
$e_{p}=\left[(e_{p}^{1})\t,\cdots,(e_{p}^{N})\t\right]\t$, $e_{r}^{i}=[(e_{\varpi}^{i})\t,(e_{p}^{i})\t]\t$
and $e_{r}=\left[(e_{\varpi})\t,(e_{p})\t\right]\t$. Note that $\xi^{i}=[(x^{i})\t,(y^{i})\t]\t=[0,(z^{*})\t]\t+[(x^{i})\t,(e_{y}^{i})\t]\t+[0,(e_{\varpi}^{i})\t]\t$
and $\zeta^{i}=\zeta^{*,i}+e_{r}^{i}$ . Then the closed-loop system composed
of (\ref{eq:nominal_sys}), (\ref{eq:zeta}), and (\ref{eq:con})
can be cast into the error dynamics as follows
\begin{gather}
\dot{e}_{r}  =\tilde{h}_{\zeta}\left(e_{r},\mu\right),\quad \dot{e}_{s}^{i}  =g_{c}^{i}\left(e_{s}^{i},e_{r},d^{i},\mu\right), \label{eq:dot_e_r}\\
e^{i}  =[0,I,0]e_{s}^{i}+[I,0]e_{r}^{i},\quad u^{i}  =\tilde{h}_{u}^{i}(e_{s}^{i},e_{r},\mu),\;i\in\mathcal{V}.\label{eq:e_i-1}
\end{gather}
where $\tilde{h}_{\zeta}$ and $g_{c}^{i}$ in (\ref{eq:dot_e_r})
can be derived from the definition, and
$ g_{c}^{i}(e_{s}^{i},e_{r},d^{i},\mu)
 = \big[\tilde{g}_{x}^{i}(e_{s}^{i},e_{r},d^{i},\mu)\t, \tilde{g}_{y}^{i}(e_{s}^{i},e_{r},d^{i},\mu)\t, \tilde{h}_{\varsigma}^{i}(e_{s}^{i},e_{r},\mu)\t \big]\t$. 
 
\begin{figure}
\includegraphics[viewport=10bp 0bp 800bp 200bp,clip,scale=0.5]{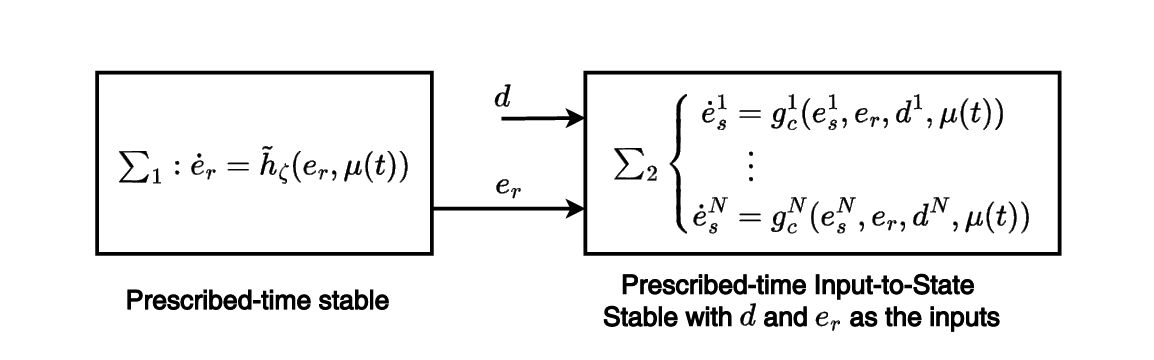}\caption{\label{fig:CasSys} Cascaded system $\Sigma=[\Sigma_{1}\t,\Sigma_{2}\t]\t$
with $d=\left[(d^{1})\t,\cdots,(d^{N})\t\right]\t$.}
\end{figure}

As illustrated in Fig. \ref{fig:CasSys},
the error system is in a cascaded form. With the decomposition of
$e^{i}$ in (\ref{eq:e_i-1}), in order to show (\ref{eq:objective}),
it suffices to prove the prescribed-time stability of $e_{r}$- and
$e_{s}^{i}$-dynamics, i.e., there exist $\mathcal{KL}_{T}^{e}$ functions
$\beta_{r}$, $\beta_{s}^{i}$ such that
\begin{equation}
\begin{aligned}\|e_{r}(t)\| & \leq\beta_{r}(\|e_{r}(t_{0})\|,t-t_{0}),\\
\|e_{s}^{i}(t)\| & \leq\beta_{s}^{i}(\|e_{s}^{i}(t_{0})\|,t-t_{0}),\quad i\in\mathcal{V}.
\end{aligned}
\label{eq:objective-1}
\end{equation}

\subsection{Prescribed-time Stabilization of Cascaded System }

\subsubsection{Changing Lyapunov Function Method}

We propose three conditions sufficient for prescribed-time stabilization
of the cascaded system (\ref{eq:dot_e_r})-(\ref{eq:e_i-1}).

$\textbf{C}_{\textbf{1}}$: The $e_{r}$-dynamics in (\ref{eq:dot_e_r})
admits a prescribed-time Lyapunov function $V_{r}(e_{r}):\mathbb{R}^{mN}\to\mathbb{R}_{\geq0}$
such that
$
V_{r}(e_{r})\sim\{\underline{\alpha}_{r},\bar{\alpha}_{r},\tilde{\alpha}_{r}|\dot{e}_{r}=\tilde{h}_{\zeta}\left(e_{r},\mu\right)\}
$
holds;

$\textbf{C}_{\textbf{2}}$: The $e_{s}^{i}$-dynamics in (\ref{eq:dot_e_r})
admits a prescribed-time ISS Lyapunov function $V_{s}^{i}(e_{s}^{i}):\mathbb{R}^{s}\to\mathbb{R}_{\geq0}$
such that $V_{s}^{i}(\text{\ensuremath{e_{s}^{i}}})\sim\{ \underline{\alpha}_{s}^{i},\bar{\alpha}_{s}^{i},\tilde{\alpha}_{s}^{i},[\sigma_{r}^{i},\tilde{\sigma}_{r}^{i}],[\sigma_{d}^{i},\tilde{\sigma}_{d}^{i}]|\dot{e}_{s}^{i}=g_{c}^{i}(e_{s}^{i},e_{r},d^{i},\mu)\},i \in\mathcal{V}$
holds for some $\sigma_{d}^{i}\in\mathcal{K}_{\infty}^{e}$ and $\sigma_{r}^{i}\in\mathcal{K}_{\infty}$;

$\textbf{C}_{\textbf{3}}$: The function $\tilde{h}_{\zeta}$ in (\ref{eq:dot_e_r})
satisfies $\Vert\tilde{h}_{\zeta}\left(e_{r},\mu\right)\Vert\leq\gamma_{\zeta}(\mu)\Vert e_{r}\Vert$
for some $\gamma_{\zeta}\in\mathcal{K}_{\infty}$; $\tilde{h}_{\varsigma}^{i}$
in $g^i_c$ and $\tilde{h}_{u}^{i}$ in (\ref{eq:e_i-1}) satisfy $\Vert\tilde{h}_{\varsigma}^{i}(e_{s}^{i},e_{r},\mu)\Vert\leq\gamma_{\varsigma}^{i}(\mu)\Vert[(e_{s}^{i})\t,e_{r}\t]\Vert$ and $\Vert\tilde{h}_{u}^{i}(e_{s}^{i},e_{r},\mu)\Vert\leq\gamma_{u}^{i}(\mu)\Vert[(e_{s}^{i})\t,e_{r}\t]\Vert$ for $i\in\mathcal{V}$ 
and  some $\gamma_{\varsigma}^{i},\gamma_{u}^{i}\in\mathcal{K_{\infty}}$.

Note that condition $\textbf{C}_{\textbf{1}}$ implies that $\dot{V}_{r}\leq-\tilde{\alpha}(\mu)V_{r}$.
Invoking comparison lemma leads to $V_{r}(e_{r}(t))\leq V_{r}(e_{r}(t_{0}))\kappa^{-1}(\tilde{\alpha}_{r}(\mu))$
where $\kappa^{-1}(\tilde{\alpha}(\mu))$ is denoted in (\ref{eq:kappa_t}).
Due to $V_{r}(e_{r})\geq\underline{\alpha}_{r}(\|e_{r}\|)$, it gives
\begin{equation}
\Vert e_{r}(t)\Vert\leq\underline{\alpha}_{r}^{-1}\left(\bar{\alpha}_{r}(\Vert e_{r}(t_{0})\Vert)\kappa^{-1}(\tilde{\alpha}_{r}(\mu))\right),\label{eq:beta_r-2}
\end{equation}
showing that the state of the first subsystem goes to zero at the prescribed time
$t_{0}+T$ and the first inequality in (\ref{eq:objective-1}) is
achieved. In order to investigate how the $e_{r}$-dynamics affects
the convergence of $e_{s}^{i}$-dynamics, we introduce the change
of the Lyapunov function for the $e_{s}^{i}$-dynamics as
\begin{equation}
W^{i}(\mu,e_{s}^{i})=\sigma_{s}^{i}(\mu)V_{s}^{i}(e_{s}^{i})\label{eq:W^i}
\end{equation}
with $\sigma_{s}^{i}\in\mathcal{K}_{\infty}$. Then, the prescribed-time convergence
result for the whole system is given in the following theorem.
\bthm \label{thm:1} Consider the system composed of (\ref{eq:nominal_sys}),
(\ref{eq:zeta}) and (\ref{eq:con}). Suppose the
closed-loop system (\ref{eq:dot_e_r})-(\ref{eq:e_i-1}) after the
state transformation satisfies conditions\textbf{ $\textbf{C}_{\textbf{1}}$}-$\textbf{C}_{\textbf{3}}$.
Define functions $\varrho_{r}(c,s)=\underline{\alpha}_{r}^{-1}\left(c\exp\left(-\sint_{0}^{s}\tau^{-2}\tilde{\alpha}_{r}(\tau)\mathrm{d}\tau\right)\right)$
and $\varrho_{s}^{i}(c,s)=\underline{\alpha}_{s}^{i,-1}(\gamma_{s}^{e,i}(c)/(\sigma_{s}^{i}(s)))$
with some $\gamma_{s}^{e,i}\in\mathcal{K}_{\infty}^{e}$ and $c\geq0$.
Suppose
\begin{equation}
\max\{\gamma_{\zeta}(s),\gamma_{\varsigma}^{i}(s),\gamma_{u}^{i}(\mu)\}=\mathcal{S}[1/(\varrho_{r}(c,s))]\label{eq:con1}
\end{equation}
and there exists a $\mathcal{K}_{\infty}$ function $\sigma_{s}^{i}$
for (\ref{eq:W^i}) such that
\begin{gather}
\frac{\mathrm{d}\sigma_{s}^{i}(s)}{\mathrm{d}s}\leq s^{-2}\sigma_{s}^{i}(s)\tilde{\alpha}_{s}^{i}(s)/2,\label{eq:con2}\\
\max\{\gamma_{\varsigma}^{i}(s),\gamma_{u}^{i}(s)\}=\mathcal{S}[1/(\varrho_{s}^{i}(c,s))],\label{eq:S1}\\
\sigma_{s}^{i}(s)\tilde{\sigma}_{d}^{i}(s)=\mathcal{S}[\exp(c)\tilde{\alpha}_{s}^{i}(s)],\label{eq:S2}\\
\sigma_{s}^{i}(s)\tilde{\sigma}_{r}^{i}(s)=\mathcal{S}[\tilde{\alpha}_{s}^{i}(s)/\sigma_{r}^{i}(\varrho_{r}(c,s))]\label{eq:con3}
\end{gather}
hold. Then, the problem of DPTCO  is solved for any bounded initial condition. \ethm
\proofnow Due to $\textbf{C}_{\textbf{2}}$, one has $\dot{V_{s}^{i}}(e_{s}^{i})\leq-\tilde{\alpha}_{s}^{i}(\mu)V_{s}^{i}(e_{s}^{i})+\tilde{\sigma}_{d}^{i}(\mu)\sigma_{d}^{i}(\Vert d^{i}\Vert)+\tilde{\sigma}_{r}^{i}(\mu)\sigma_{r}^{i}(\Vert e_{r}\Vert)$.
Taking time derivative of $W^{i}(\mu,e_{s}^{i})$ in (\ref{eq:W^i})
and using (\ref{eq:con2}) yields $\dot{W}^{i}\leq-\tilde{\alpha}_{s}^{i}(\mu)W^{i}/2+\tilde{\gamma}_{d}^{i}(\mu)\sigma_{d}^{i}(\Vert d^{i}\Vert)+\tilde{\gamma}_{r}^{i}(\mu)\sigma_{r}^{i}(\Vert e_{r}\Vert) $, 
where $\tilde{\gamma}_{d}^{i}(\mu)=\sigma_{s}^{i}(\mu)\tilde{\sigma}_{d}^{i}(\mu)$
and $\tilde{\gamma}_{r}^{i}(\mu)=\sigma_{s}^{i}(\mu)\tilde{\sigma}_{r}^{i}(\mu)$.
Invoking comparison lemma yields
\begin{gather}
W^{i}(t)\leq W^{i}(t_{0})\kappa^{-\frac{1}{2}}\left(\tilde{\alpha}_{s}^{i}(\mu)\right)\nonumber \\
+\kappa^{-\frac{1}{2}}(\tilde \alpha_s^i (\mu))\int_{t_{0}}^{t}\kappa^{\frac{1}{2}}(\tilde \alpha_s^i (\mu(\tau )))\tilde{\gamma}_{d}^{i}(\mu(\tau))\sigma_{d}^{i}(\Vert d^{i}(\tau)\Vert)\mathrm{d}\tau\nonumber \\
+\kappa^{-\frac{1}{2}}(\tilde \alpha_s^i (\mu))\int_{t_{0}}^{t}\kappa^{\frac{1}{2}}(\tilde \alpha_s^i (\mu(\tau )))\tilde{\gamma}_{r}^{i}(\mu(\tau))\sigma_{r}^{i}(\Vert e_{r}(\tau)\Vert)\mathrm{d}\tau.\notag
\end{gather}
Denote the bound of $\|d^{i}\|$ as $\bar{d}^{i}$. Given (\ref{eq:S2})
with $c=0$, one has $\sup_{\mu\in R_{p}}\left(\tilde{\gamma}_{d}^{i}(\mu)/\tilde{\alpha}_{s}^{i}(\mu)\right)<\infty$.
As a result, the second term in $W^i(t)$
can be calculated as
\begin{align}
&	\kappa^{-\frac{1}{2}}(\tilde \alpha_s^i (\mu))\int_{t_{0}}^{t}\kappa^{\frac{1}{2}}(\tilde \alpha_s^i (\mu(\tau )))\tilde{\gamma}_{d}^{i}(\mu(\tau))\sigma_{d}^{i}(\Vert d^{i}(\tau)\Vert)\mathrm{d}\tau\notag \\
& \leq \sigma_{d}^{i}(\bar{d}^{i})\kappa^{-\frac{1}{2}}\left(\tilde{\alpha}_{s}^{i}(\mu)\right)\int_{t_{0}}^{t}\kappa^{\frac{1}{2}}(\tilde \alpha_s^i (\mu(\tau )))\tilde{\gamma}_{d}^{i}(\mu(\tau))\mathrm{d}\tau \notag \\
& = \sigma_{d}^{i}(\bar{d}^{i})\kappa^{-\frac{1}{2}}\left(\tilde{\alpha}_{s}^{i}(\mu)\right)\int_{t_{0}}^{t}\kappa^{\frac{1}{2}}(\tilde \alpha_s^i (\mu(\tau ))) \notag\\
&\quad \times \left(\tilde{\gamma}_{d}^{i}(\mu(\tau))/\tilde \alpha_s^i (\mu(\tau ))\right)\mathrm{d}\left(\sint_{t_{0}}^{\tau}\tilde{\alpha}_{s}^{i}(\mu(s))\mathrm{d}s\right) \notag \\
 & \leq \epsilon_{d}^{i}\big(1-\kappa^{-\frac{1}{2}}\left(\tilde{\alpha}_{s}^{i}(\mu)\right)\big),\label{eq:SecTerm-1}
\end{align}
where $\epsilon_{d}^{i}=2\sigma_{d}^{i}(\bar{d}^{i})\sup_{\mu\in\mathbb{R}_{p}}\left(\tilde{\gamma}_{d}^{i}(\mu)/\tilde{\alpha}_{s}^{i}(\mu)\right)$
is a finite constant, and we used the facts that $\int_{t_0}^t\tilde \alpha^i_s (\mu(\tau))\mathrm d\tau =\int_{t_0}^t \mathrm d \big(\sint_{t_0}^\tau \tilde \alpha^i_s (\mu(s))\mathrm ds\big)$ and $\int_{t_{0}}^{t}\kappa^{\frac{1}{2}}(\tilde{\alpha}_{s}^{i}(\mu(\tau)))\mathrm{d}\big(\sint_{t_{0}}^{\tau}\tilde{\alpha}_{s}^{i}(\mu(s))\mathrm{d}s\big) =2 \big(\kappa^{\frac{1}{2}}(\tilde{\alpha}_{s}^{i}(\mu))-1\big)$. 
By (\ref{eq:beta_r-2}), one has
\begin{align}
\|e_{r}(t)\| & 
  \leq \underline{\alpha}_{r}^{-1}\left(\epsilon\bar{\alpha}_{r}(\Vert e_{r}(t_{0})\Vert)\exp\left(-\sint_{0}^{\mu(t)}\tau^{-2}\tilde{\alpha}_{r}(\tau)\mathrm{d}\tau\right)\right)\nonumber \\
 & =\varrho_{r}(\epsilon\bar{\alpha}_{r}(\Vert e_{r}(t_{0})\Vert),\mu(t)),\label{eq:varrho_r}
\end{align}
where $\epsilon=\exp\left(\sint_{0}^{\mu(t_{0})}\tau^{-2}\tilde{\alpha}_{r}(\tau)\mathrm{d}\tau\right)$.
Similar to (\ref{eq:SecTerm-1}), due to (\ref{eq:con3}) and (\ref{eq:varrho_r}),
the third term in $W^i(t)$ satisfies
\begin{align}
	& \kappa^{-\frac{1}{2}}(\tilde \alpha_s^i (\mu))\int_{t_{0}}^{t}\kappa^{\frac{1}{2}}(\tilde \alpha_s^i (\mu(\tau )))\tilde{\gamma}_{r}^{i}(\mu(\tau))\sigma_{r}^{i}(\Vert e_{r}(\tau)\Vert)\mathrm{d}\tau\nonumber \\
	& \leq\epsilon_{r}^{i}\big(1-\kappa^{-\frac{1}{2}}\left(\tilde{\alpha}_{s}^{i}(\mu)\right)\big),\notag 
\end{align}
where $\epsilon_{r}^{i}=2\sup_{\mu\in\mathbb{R}_{p}}(\tilde{\gamma}_{r}^{i}(\mu)\sigma_{r}^{i}(\varrho_{r}(\epsilon\bar{\alpha}_{r}(\Vert e_{r}(t_{0})\Vert),\mu))/\tilde{\alpha}_{s}^{i}(\mu))$
is a finite constant. Consequently, $	W^{i}(t) \leq W^{i}(t_{0})+\epsilon_{d}^{i}+\epsilon_{r}^{i}$. 
Then according to (\ref{eq:W^i}), $e_{s}^{i}$ satisfies
\begin{align}
\Vert e_{s}^{i}(t)\Vert \leq \varrho_{s}^{i}(\|e_{s}^{i}(t_{0}\|,\mu),\label{eq:converge_e_s}
\end{align}
where $\gamma_{s}^{e,i}(\|e_{s}^{i}(t_{0})\|)=\sigma_{s}^{i}(\mu(t_{0}))\bar{\alpha}_{s}^{i}(\|e_{s}^{i}(t_{0})\|)+\epsilon_{d}^{i}+\epsilon_{r}^{i}$. Inequality  (\ref{eq:converge_e_s}) means that the
second equation in (\ref{eq:objective-1}) is achieved. As a result,
the DPTCO is achieved.

Next, we prove the boundedness of $h_{\zeta}^{i}$, $h_{\varsigma}^{i}$, and 
$h_{u}^{i}$. By (\ref{eq:beta_r-2}), (\ref{eq:con1}), (\ref{eq:S1})
and (\ref{eq:converge_e_s}), $\max\{\gamma_{\zeta}(\mu),\gamma_{\varsigma}^{i}(\mu),\gamma_{u}^{i}(\mu)\}\Vert e_{r}\Vert\leq\varepsilon_{r}$,
$\max\{\gamma_{\varsigma}^{i}(\mu),\gamma_{u}^{i}(\mu)\}\Vert e_{s}^{i}\Vert\leq\varepsilon_{s}$
hold for some finite constants $\varepsilon_{r}$, $\varepsilon_{s}$,
and $\mu\in\mathbb{R}_{p}$. Since $\tilde{h}_{\zeta}^{i}$, $\tilde{h}_{\varsigma}^{i}$, and
$\tilde{h}_{u}^{i}$ satisfy $\textbf{C}_{\textbf{3}}$,
these inequalities imply that $h_{\zeta}^{i}$, $h_{\varsigma}^{i}$, and 
$h_{u}^{i}$ are bounded for $t\in\mathcal{T}_{p}$. This completes
the proof. \eproof

\subsubsection{Time-varying State Transformation \label{subsec:TST}}

A common practice in the literature of prescribed-time control \cite{zuo2022adaptive35,song2017time36}
is the time-varying state transformation technique. When \textbf{$\textbf{C}_{\textbf{2}}$}
is not feasible, we can seek a time-varying state transformation
\begin{equation}
\tilde{e}_{s}^{i}=h_{s}^{i}(e_{s}^{i},\mu),\label{eq:til_e_s}
\end{equation}
where $h_{s}^{i}:\mathbb{R}^{n+m+m_{\varsigma}}\times\mathbb{R}_{p}\to\mathbb{R}^{s}$
is a differentiable function. Generally, the mapping from $e_{s}^{i}$
to $\tilde{e}_{s}^{i}$ is nonlinear. The $\tilde{e}_{s}^{i}$-dynamics
becomes
\begin{align*}
\dot{\tilde{e}}_{s}^{i} & =\tilde{h}_{s}^{i}(\tilde{e}_{s}^{i},e_{r},d^{i},\mu)= \frac{\partial h_{s}^{i}}{\partial e_{s}^{i}}g_{c}^{i}\left(e_{s}^{i},e_{r},d^{i},\mu\right)+\frac{\partial h_{s}^{i}}{\partial\mu}\mu^{2},
\end{align*}
where we used $\dot{\mu}=\mu^{2}$. Due to the nonlinearity of $h_{s}^{i}(\cdot)$,
$\tilde{e}_{s}^{i}=0$ may not guarantee $\dot{\tilde{e}}_{s}^{i}=0$.
With the time-varying state transformation, the closed-loop system
composed of (\ref{eq:nominal_sys}), (\ref{eq:zeta}), and (\ref{eq:con})
can be cast into the error dynamics as follows
\begin{align}
\dot{e}_{r} & =\tilde{h}_{\zeta}\left(e_{r},\mu\right),\quad 
\dot{\tilde{e}}_{s}^{i}  =\tilde{h}_{s}^{i}(\tilde{e}_{s}^{i},e_{r},d^{i},\mu),\;i\in\mathcal{V}
\label{eq:dot_e_r-1}
\end{align}
and $\varsigma^{i}$-dynamics and $u^{i}$ in (\ref{eq:zeta}), (\ref{eq:con})
can be rewritten as
\begin{align}
\dot{\varsigma}^{i} & =\bar{h}_{\varsigma}^{i}(\tilde{e}_{s}^{i},e_{r},\mu),\quad 
u^{i}  =\bar{h}_{u}^{i}(\tilde{e}_{s}^{i},e_{r},\mu) \label{eq:bar_h_sigma}
\end{align}
with some functions $\bar{h}_{\varsigma}^{i}$ and $\bar{h}_{u}^{i}$
derived from (\ref{eq:zeta}), (\ref{eq:con}) and (\ref{eq:til_e_s}).
Similarly, $\tilde{e}_{s}^{i}=0$ may not guarantee $\bar{h}_{\varsigma}^{i}=0$
and $\bar{h}_{u}^{i}=0$. We modify conditions \textbf{$\textbf{C}_{\textbf{2}}$},
$\textbf{C}_{\textbf{3}}$ to \textbf{$\textbf{C}_{\textbf{2}}'$},
\textbf{$\textbf{C}_{\textbf{3}}'$} as follows.

\textbf{$\textbf{C}_{\textbf{2}}'$}: There exists a time-varying
state transformation (\ref{eq:til_e_s}) such that $\tilde{e}_{s}^{i}$-dynamics
in (\ref{eq:dot_e_r-1}) admits a prescribed-time ISS Lyapunov function
$\tilde{V}_{s}^{i}(\tilde{e}_{s}^{i}):\mathbb{R}^{s}\to\mathbb{R}_{\geq0}$
and $\tilde{V}_{s}^{i}(\tilde{e}_{s}^{i})\sim\{\underline{\alpha}_{\tilde{s}}^{i},\bar{\alpha}_{\tilde{s}}^{i},\tilde{\alpha}_{\tilde{s}}^{i},[\sigma_{\tilde{r}}^{i},\tilde{\sigma}_{\tilde{r}}^{i}],[\sigma_{\tilde{d}}^{i},\tilde{\sigma}_{\tilde{d}}^{i}]|\dot{\tilde{e}}_{s}^{i}=\tilde{h}_{s}^{i}(\tilde{e}_{s}^{i},e_{r},d^{i},\mu)\},\forall i\in\mathcal{V}$
holds for some $\sigma_{\tilde{r}}^{i}\in\mathcal{K}_{\infty}$ and
$\sigma_{\tilde{d}}^{i}\in\mathcal{K}_{\infty}^{e}$. \textit{Moreover,
the boundedness of $\tilde{e}_{s}^{i}$ implies prescribed-time convergence
of $e_{s}^{i}$.}

$\textbf{C}_{\textbf{3}}'$: $\bar{h}_{\varsigma}^{i}$  
and $\bar{h}_{u}^{i}$ in (\ref{eq:bar_h_sigma}) satisfy $\Vert\bar{h}_{\varsigma}^{i}(\tilde{e}_{s}^{i},e_{r},\mu)\Vert\leq\gamma_{s}^{i}(\Vert\tilde{e}_{s}^{i}\Vert_{\text{\ensuremath{\mathcal{T}}}})+\tilde{\gamma}_{\varsigma}^{i}(\mu)\Vert e_{r}\Vert$ and $\Vert\bar{h}_{u}^{i}(\tilde{e}_{s}^{i},e_{r},\mu)\Vert\leq\tilde{\gamma}_{s}^{i}(\Vert\tilde{e}_{s}^{i}\Vert_{\mathcal{T}})+\tilde{\gamma}_{u}^{i}(\mu)\Vert e_{r}\Vert,i\in\mathcal{V}$
 for $\mu\in\mathbb{R}_{p}$ where $\|\tilde{e}_{s}^{i}\|_{\mathcal{T}}=\sup_{t\in\mathcal{T}_{p}}\|\tilde{e}_{s}^{i}(t)\|$,
$\gamma_{s}^{i}$, $\tilde{\gamma}_{s}^{i}\in\mathcal{K}_{\infty}^{e}$
and $\tilde{\gamma}_{u}^{i}$, $\tilde{\gamma}_{\varsigma}^{i}\in\mathcal{K}_{\infty}$.

\bthm \label{thm:2} Consider the system composed of (\ref{eq:nominal_sys}),
(\ref{eq:zeta}), and (\ref{eq:con}). Suppose the
closed-loop system (\ref{eq:e_i-1}) and (\ref{eq:dot_e_r-1}) 
after the state transformation satisfies conditions $\text{\textbf{C}}_{\textbf{1}}$,\textbf{
$\textbf{C}_{\textbf{2}}'$}, \textbf{$\textbf{C}_{\textbf{3}}'$}
with
\begin{equation}
\max\{\gamma_{\zeta}(s),\tilde{\gamma}_{\varsigma}^{i}(s),\tilde{\gamma}_{u}^{i}(s)\}=\mathcal{S}[1/(\varrho_{r}(c,s))],\label{eq:con1-1}
\end{equation}
where $\varrho_{r}(c,s)$ is defined in Theorem \ref{thm:1}, and
\begin{gather}
\tilde{\sigma}_{\tilde{d}}^{i}(s)=\mathcal{S}[\exp(c)\tilde{\alpha}_{\tilde{s}}^{i}(s)],\label{eq:con2-1}\\
\tilde{\sigma}_{\tilde{r}}^{i}(s)=\mathcal{S}[\tilde{\alpha}_{\tilde{s}}^{i}(s)/(\sigma_{\tilde{r}}^{i}(\varrho_{r}(c,s)))]\label{eq:con3-3}
\end{gather}
hold. Then, the problem of DPTCO  is solved for any bounded initial condition. \ethm

\proofnow Due to \textbf{$\textbf{C}_{\textbf{2}}'$}, one has $\dot{\tilde{V}}_{s}^{i}(\tilde{e}_{s}^{i})\leq-\tilde{\alpha}_{\tilde{s}}^{i}(\mu)\tilde{V}_{s}^{i}(\tilde{e}_{s}^{i})+\tilde{\sigma}_{\tilde{d}}^{i}(\mu)\sigma_{\tilde{d}}^{i}(\Vert d^{i}\Vert)+\tilde{\sigma}_{\tilde{r}}^{i}(\mu)\sigma_{\tilde{r}}^{i}(\Vert e_{r}\Vert)$.
Invoking comparison lemma yields
\begin{gather*}
\tilde{V}_{s}^{i}(\tilde{e}_{s}^{i}(t))\leq\tilde{V}_{s}^{i}(\tilde{e}_{s}^{i}(t_{0}))\kappa^{-1}(\tilde{\alpha}_{\tilde{s}}^{i}(\mu))\\
+\int_{t_{0}}^{t}\exp\left(-\sint_{\tau}^{t}\tilde{\alpha}_{\tilde{s}}^{i}(\mu(s))\mathrm{d}s\right)\tilde{\sigma}_{\tilde{d}}^{i}(\mu(\tau))\sigma_{\tilde{d}}^{i}(\Vert d^{i}(\tau)\Vert)\mathrm{d}\tau\\
+\int_{t_{0}}^{t}\exp\left(-\sint_{\tau}^{t}\tilde{\alpha}_{\tilde{s}}^{i}(\mu(s))\mathrm{d}s\right)\tilde{\sigma}_{\tilde{r}}^{i}(\mu(\tau))\sigma_{\tilde{r}}^{i}(\Vert e_{r}(\tau)\Vert)\mathrm{d}\tau.
\end{gather*}
Similar to the deviations in (\ref{eq:SecTerm-1}),
by (\ref{eq:con2-1}) and (\ref{eq:con3-3}) , the bound of $\tilde{V}_{s}^{i}$
satisfies $\tilde{V}_{s}^{i}(\tilde{e}_{s}^{i}(t))\leq\tilde{V}_{s}^{i}(\tilde{e}_{s}^{i}(t_{0}))\kappa^{-1}(\tilde{\alpha}_{\tilde{s}}^{i}(\mu)) +(\tilde{\epsilon}_{d}^{i}+\tilde{\epsilon}_{r}^{i})(1-\kappa^{-1}(\tilde{\alpha}_{\tilde{s}}^{i}(\mu)))$,
where $\tilde{\epsilon}_{d}^{i}=\sup_{\mu\in\mathbb{R}_{p}}\big(\tilde{\sigma}_{\tilde{d}}^{i}(\mu)/\tilde{\alpha}_{\tilde{s}}^{i}(\mu)\big)\sigma_{\tilde{d}}^{i}(\bar{d}^{i})<\infty$ and $\tilde{\epsilon}_{r}^{i}=\sup_{\mu\in\mathbb{R}_{p}}(\tilde{\sigma}_{\tilde{r}}^{i}(\mu)\sigma_{\tilde{r}}^{i}(\varrho_{r}(\epsilon\bar{\alpha}_{r}(\Vert e_{r}(t_{0})\Vert),\mu))/\tilde{\alpha}_{\tilde{s}}^{i}(s))<\infty$.

The inequality  implies that $\tilde{e}_{s}^{i}$
is bounded. Since the boundedness of $\tilde{e}_{s}^{i}$ implies
the prescribed-time convergence of $e_{s}^{i}$ by condition $\textbf{C}_{\textbf{2}}'$,
the second equation in (\ref{eq:objective-1}) is achieved and outputs
of the agents converge to the optimum within prescribed time.
Similar to the proof of Theorem \ref{thm:1}, by (\ref{eq:con1-1}),
we have $\max\{\gamma_{\zeta}(\mu),\tilde{\gamma}_{\varsigma}^{i}(\mu),\tilde{\gamma}_{u}^{i}(\mu)\}\Vert e_{r}\Vert<\infty$,
and then the boundedness of all signals is guaranteed. \eproof

\brem 
The most fundamental difference between the changing Lyapunov function method and the time-varying state transformation method lies in their distinct approaches to deriving the prescribed-time convergence of $e_s^i$. The changing Lyapunov function method seeks a transformation of the Lyapunov function $V_s^i (e_s^i)$ in $\textbf{C}_{\textbf{2}}$, as exemplified  by $W^i (\mu, e_s^i)$ in \eqref{eq:W^i}. By showing the boundedness of  $W^i (\mu, e_s^i)$, one obtains the prescribed-time convergence of  $e_s^i$. This is a process that leverages the common idea that proving boundedness of a Lyapunov function is more feasible than directly proving its convergence. However, when $e_s^i$-dynamics is nonlinear and has a relative degree greater than one, it often becomes necessary to employ some control techniques such as sliding-mode control or backstepping method, making it very challenging for $e_s^i$-dynamics in \eqref{eq:dot_e_r} to admit a prescribed-time ISS Lyapunov function $V^i_s(e^i_s)$ in $\textbf{C}_{\textbf{2}}$. In such cases, one can instead apply the time-varying state transformation method. 
\erem 

\brem
The core of time-varying state transformation method  lies in identifying a nonlinear time-varying transformation \eqref{eq:til_e_s} for $e_s^i$, converting the prescribed-time convergence problem of $e_s^i$ into a boundedness problem for a new variable $\tilde e_s^i$, whose dynamics admit a prescribed-time ISS Lyapunov function $\tilde V_s^i (\tilde e_s^i)$ in \textbf{$\textbf{C}_{\textbf{2}}'$}. Although the changing Lyapunov function method is theoretically intuitive, the requirement that $e_s^i$-dynamics in \eqref{eq:dot_e_r} admits a prescribed-time ISS Lyapunov function in $\textbf{C}_{\textbf{2}}$ is often difficult to satisfy.  In contrast, the time-varying state transformation method, though more technically complex, particularly in designing the nonlinear time-varying transformation \eqref{eq:til_e_s}, is more generally applicable to higher-order nonlinear systems. 
\erem

\section{Prescribed-time Optimum Seeking \label{sec:Optimal-Trajectory-Generator}}
In this section, we elaborate on the design of $\zeta^{i}$-dynamics.
The two subsystems of $\zeta^{i}$-dynamics, namely $\varpi^{i}$-
and $p^{i}$-dynamics, are designed as,
\begin{align}
\dot{\varpi}^{i} & =-\alpha(\mu)\left(\ssum_{j\in\mathcal{N}^{i}}(\varpi^{i}-\varpi^{j})+\nabla f^{i}(\varpi^{i})+p^{i}\right),\label{eq:dot_varpi_i-1}\\
\dot{p}^{i} & =\alpha(\mu)\ssum_{j\in\mathcal{N}^{i}}(\varpi^{i}-\varpi^{j}),\;i\in\mathcal{V},\label{eq:dot_p_i}
\end{align}
where $\alpha\in\mathcal{K}_{\infty}$ is a differentiable function
to be designed.

Let $r=1_{N}/\sqrt{N}\in\mathbb{R}^{N}$ and $R\in\mathbb{R}^{N\times(N-1)}$
be such that $r\t R=0$, $R\t R=I_{N-1}$. Therefore, $RR\t=\Pi_{N}=I_{N}-\frac{1}{N}1_{N}1_{N}\t$
and $[r,R]$ is an orthogonal matrix. Define $\mathcal{L}_{R}=R\t\mathcal{L}R$,
$\bar{\mathcal{L}}=\mathcal{L}\otimes I_{m}$, $\bar{\mathcal{L}}_{R}=\mathcal{L}_{R}\otimes I_{m}$,
$\bar{R}=R\otimes I_{m}$ and $\bar{r}=r\otimes I_{m}$.For a connected graph, $\mathcal{L}_{R}$ is a positive matrix and
$\lambda_{2}I_{N-1}\leq\mathcal{L}_{R}\leq\lambda_{N}I_{N-1}$ where
$\lambda_{2}$ and $\lambda_{N}$ are the second smallest and largest
eigenvalues of $\mathcal{L}$, respectively. Let $\varpi=[(\varpi^{1})\t,\cdots,(\omega^{N})\t]\t$
and $p=[(p^{1})\t,\cdots,(p^{N})\t]\t$. The dynamics (\ref{eq:dot_varpi_i-1})
and (\ref{eq:dot_p_i}) for the group of agents can be written compactly
as
\begin{align}
\dot{\varpi} & =-\alpha(\mu)\bar{\mathcal{L}}\varpi-\alpha(\mu)\nabla F(\varpi)-\alpha(\mu)p,\label{eq:dot_varpi}\\
\dot{p} & =\alpha(\mu)\bar{\mathcal{L}}\varpi,\label{eq:dot_p}
\end{align}
 where $\nabla F(\varpi)=[\nabla f^{1}(\varpi^{1});\cdots;\nabla f^{N}(\varpi^{N})]$.
Note that the systems (\ref{eq:dot_varpi}) and (\ref{eq:dot_p}) are
in the form of (\ref{eq:zeta}).

\bproposition \label{prop:sol} Consider (\ref{eq:dot_varpi}) and
(\ref{eq:dot_p}) under Assumptions \ref{ass:graph}, \ref{ass:cost_func}, 
and \ref{ass:solvable}. Let $z^{*}$ satisfy \eqref{eq:y_opt} and thus $1_{N}\otimes z^{*}$
is the optimum of  the optimization problem (\ref{eq:opti_problem}).
Then
\begin{align}
\varpi^{*} & =1_{N}\otimes z^{*},\quad 
p^{*}  =-\nabla F(1_{N}\otimes z^{*})\label{eq:varpi*}
\end{align}
is the solution of
\begin{equation}
0=-\bar{\mathcal{L}}\varpi-\nabla F(\varpi)-p,\quad0=\bar{\mathcal{L}}\varpi\label{eq:equi-varpi}
\end{equation}
when the  initial value of $p^{i}(t_{0})$ satisfies $\ssum_{i=1}^{N}p^{i}(t_{0})=0$. 
\eproposition

\proofnow
From (\ref{eq:equi-varpi}), the solution
satisfies
\begin{equation}
	\varpi^{*}\in\mbox{span}\{1_{N}\otimes\upsilon\},\quad p^{*}=-\nabla F(\varpi^{*})\label{eq:varpi_p*}
\end{equation}
for any vector $\upsilon\in\mathbb{R}^{m}$. Since $\mathcal{G}$
is undirected and connected, the Laplacian matrix matrix $\mathcal{L}$
is symmetric and its null space is spanned by $1_{N}$, then $1_{N}\t\mathcal{L}=0$.
By $\ssum_{i=1}^{N}p^{i}(t_{0})=0$, left multiplying (\ref{eq:dot_p})
by $1_{N}\t\otimes I_{m}$ yields
\begin{gather}
	\ssum_{i=1}^{N}\dot{p}^{i}(t)=0 \notag \\
	\implies\ssum_{i=1}^{N}p^{i}(t)=\ssum_{i=1}^{N}p^{i}(t_{0})=0,\;t\geq t_{0}. \label{eq:initial_p}
\end{gather} 
Then left multiplying the first equation in (\ref{eq:equi-varpi})
by $1_{N}\t\otimes I_{m}$ yields
$
\ssum_{i=1}^{N}\nabla f^{i}(v)=0$. 
For the optimization problem (\ref{eq:opti_problem}), the necessary
and sufficient condition for a point $z^{*}$ to be the unique optimum is $\nabla f(z^{*})=\ssum_{i=1}^{N}f^{i}(z^{*})=0$, and
thus we have $v=z^{*}$ and (\ref{eq:varpi*}). Substituting (\ref{eq:varpi*})
into the second equation of (\ref{eq:varpi_p*}) leads to (\ref{eq:varpi*}).
\eproof

As introduced in Section \ref{sec:A-Cascade-Design}, we use the coordinate
transformation, i.e., $e_{\varpi}=\varpi-1_{N}\otimes z^{*}$ and $e_{p}=p+\nabla F(1_{N}\otimes z^{*})$
with $e_{\varpi}$ and $e_{p}$ being the error variables for distributed
optimal value seeking problem. From Proposition \ref{prop:sol}, (\ref{eq:dot_varpi})
and (\ref{eq:dot_p}), the $e_{r}$-dynamics can be obtained, with $e_{r}=[e_{\varpi}\t,e_{p}\t]\t$,
as
\begin{align}
\dot{e}_{\varpi} & =-\alpha(\mu)(\bar{\mathcal{L}}e_{\varpi}+\nabla\tilde{F}(e_{\varpi})+e_{p}),\label{eq:dot_e_varpi}\\
\dot{e}_{p} & =\alpha(\mu)\bar{\mathcal{L}}e_{\varpi},\label{eq:dot_e_p}
\end{align}
 where 
\begin{align*}
\nabla\tilde{F}(e_{\varpi})= & \left[(\nabla f^{1}(\varpi^{1})-\nabla f^{1}(z^{*}))\t,\right.\\
 & \left.\cdots,(\nabla f^{N}(\varpi^{N})-\nabla f^{N}(z^{*}))\t\right]\t.
\end{align*}

\bthm\label{the:auxi_genera} Consider the $\zeta^{i}$-dynamics in (\ref{eq:dot_varpi_i-1})
and (\ref{eq:dot_p_i}) under Assumptions \ref{ass:graph}, \ref{ass:cost_func}
and \ref{ass:solvable}. Define
\begin{equation}
\begin{gathered}
c_{1}  =\max\left\{ 1/\lambda_{2}, (1+2\varrho_{c}^{2})/(2\rho_{c})\right\},\\
c_{2}=c_1\min\left\{ 1/2,1/(2\lambda_{N})\right\},\\
c_{3}  =c_{1}\max\left\{ 1,1/\lambda_{2}\right\} +1,\quad c^{*}=1/4c_{3},
\end{gathered}\label{eq:c_1}
\end{equation}
where $\rho_{c}$ and $\varrho_{c}$ are given in Assumption \ref{ass:cost_func}.
If $\ssum_{i=1}^{N}p^{i}(t_{0})=0$ and
\begin{equation}
\frac{\mathrm{d}\alpha(s)}{\mathrm{d}s}\leq c^{*}s^{-2}(\alpha(s))^{2}/2 \label{eq:par_alpha}
\end{equation}
holds for $s\in\mathbb{R}_{\geq0}$, then the $e_{r}$-dynamics satisfies
condition $\textbf{C}_{\textbf{1}}$ with
\begin{equation}
\begin{gathered}\underline{\alpha}_{r}(s)=c_{2}s^{2},\quad\bar{\alpha}_{r}(s)=c_{3}s^{2},\\
\tilde{\alpha}_{r}(s)=2c^{*}\alpha(s),\quad\gamma_{\zeta}(s)=\max\{2\lambda_{N}+\varrho_{c},1\}\alpha(s).
\end{gathered}
\label{eq:alpha_r}
\end{equation}
 Moreover, the bounds of $e_{r}$ and $\dot{e}_{r}$ satisfy
\begin{align}
\Vert e_{r}\Vert & \leq\gamma_{r}(\|e_{r}(t_{0})\|)\kappa^{-c^{*}}\left(\alpha(\mu)\right),\label{eq:e_r_bound}\\
\|\dot{e}_{r}\| & \leq\gamma_{e}(\|e_{r}(t_{0})\|)\kappa^{-\frac{c^{*}}{2}}\left(\alpha(\mu)\right)\label{eq:er_bound-1}
\end{align}
for some $\gamma_{r},\gamma_{e}\in\mathcal{K}_{\infty}$. 
The proof is given in appendix. 
\ethm

\brem \label{rem:4-2} According
to the proof, (\ref{eq:par_alpha}) is used to guarantee the prescribed-time
convergence of $\dot{e}_{r}$. When we select $\alpha(s)=ks$ with
$k\geq2/c^{*}$, 
or $\alpha(s)=k_{1}s\exp(k_{2}s)$ with $k_{1}\geq2/c^{*}$, $k_{2}\geq0$,
$\gamma_{r}(\|e_{r}(t_{0})\|)\kappa^{-c^{*}}(\alpha(\mu))$ is a $\mathcal{KL}_{T}^{e}$
function and the first inequality in (\ref{eq:objective-1}) holds.
Since the higher-order derivative of $\varpi^{i}$ may not exist, the
prescribed-time convergence of $\dot{e}_{r}$ will play a very important
role in the design of the local prescribed-time tracking controller
later.\erem 

\section{Robust DPTCO for Chain-Integrator MASs \label{sec:DPTCO-for-Chain}}
In this section, we apply the DPTCO framework proposed in Section \ref{sec:A-Cascade-Design} to solve the robust DPTCO for a class of nonlinear MASs with uncertainties, called chain-integrator MASs of a relative degree greater than one. 

Since we deal with
the optimal tracking problem for each subsystem separately,\textit{ we
omit the superscript $i$ for simplicity when no confusion arises.} Therefore, the $i$-th subsystem is  expressed as 
\begin{equation}
\begin{aligned}\dot{x}_{q} & =x_{q+1},\;q=1,\cdots,m-1,\\
\dot{x}_{m} & =u+\varphi(x,d),\quad y  =x_{1},
\end{aligned}
\label{eq:chain_sys}
\end{equation}
where $x=[x_{1}\t,\cdots,x_{m}\t]\t\in\mathbb{R}^{mn}$ is the system
state with $x_{q}\in\mathcal{\mathbb{R}}^{n}$, $u\in\mathbb{R}^{n}$
control input, $y\in\mathbb{R}^{n}$ system output, and $d\in\mathbb{D}$
the uncertainties belonging to a compact set $\mathbb{D}\in\mathbb{R}^{n_{d}}$.
The function $\varphi:\mathbb{R}^{mn}\times\mathbb{D}\to\mathbb{R}^{n}$
is sufficiently smooth, and for each fixed $d$, it is bounded for all
$x\in\mathbb{R}^{mn}$ \cite{song2015tracking48}. 
According to \cite[Lemma 11.1]{chen2015stabilization22}, the function
$\varphi$ satisfies
\begin{equation}
\Vert\varphi(x,d)\Vert\leq h(\|d\|)\psi(x),\label{eq:psi^i}
\end{equation}
where $h\in\mathcal{K}_{\infty}$ is an unknown positive function, 
and $\psi(x)$ is a known positive function that is bounded for all
$x\in\mathbb{R}^{mn}$.Note that (\ref{eq:chain_sys})
is in the form of (\ref{eq:nominal_sys}).

We follow the framework developed in Section \ref{sec:A-Cascade-Design}
to solve the DPTCO problem. First, define the error as in (\ref{eq:e_y_varpi}),
i.e.,
\begin{equation}
e_{s}=\big[(x_{1}-\varpi^{i})\t,(x_{2})\t,(x_{m})\t\big]\t, \label{eq:e_s}
\end{equation}
where $\varpi^{i}$ is given in (\ref{eq:dot_varpi_i-1}) and $\varsigma$
is omitted in this section. Due to (\ref{eq:objective-1}) and Theorem
\ref{the:auxi_genera}, it suffices to design a controller $u$ such
that the prescribed-time stabilization is achieved for $e_{s}$.

Let $K=[k_{1},\cdots,k_{m-1}]\t\in\mathbb{R}^{m-1}$ such that
\begin{equation}
\Lambda=\left[\begin{array}{cc}
0_{m-2} & I_{m-2}\\
-k_{1} & -k_{2}\cdots-k_{m-1}
\end{array}\right]\label{eq:Lambda}
\end{equation}
 is Hurwitz, $L_{j}=j-1$ for $j=2,\cdots,m$, $\tilde{K}=[K\t,1]\t$
and $\Phi(\mu)$ is
\begin{equation}
\Phi(\mu)=\mbox{diag}\big\{ 1,(\alpha_{x}(\mu))^{-L_{2}},\cdots,(\alpha_{x}(\mu))^{-L_{m}}\big\} .\label{eq:Phimu}
\end{equation}
where $\alpha_{x}\in\mathcal{K}_{\infty}$ is a first-order differentiable
function to be designed. Since the system (\ref{eq:chain_sys}) is
nonlinear and has a relative degree greater than one, and the reference
trajectory $\varpi^{i}$ does not have the higher-order derivatives,
the traditional sliding-mode based tracking control cannot be applied
\cite{edwards1998sliding39}. Instead, we construct a new variable
$\tilde{s}$ as
\begin{align}
\tilde{s} & =k_{1}^{-1}(\tilde{K}\t\Phi(\mu)\otimes I_{n})e_{s}\nonumber \\
 & =k_{1}^{-1}\left(K\t\otimes I_{n}\right)r_{1}+k_{1}^{-1}\left(\alpha_{x}(\mu)\right)^{-L_{m}}x_{m}-\varpi^{i}\label{eq:til_s^i-2}
\end{align}
with
\begin{align}
r_{1}= & \left[x_{1}\t,(\alpha_{x}(\mu))^{-L_{2}}x_{2}\t,\cdots,(\alpha_{x}(\mu))^{-L_{m-1}}x_{m-1}\t\right]\t.\label{eq:til_r^i_1}
\end{align}
Then, we define the time-varying state transformation as
\begin{equation}
\tilde{e}_{s}=\alpha_{s}(\mu)\tilde{s},\label{eq:til_e_s^i}
\end{equation}
with a first-order differentiable $\alpha_{s}\in\mathcal{K}_{\infty}$
to be designed. By doing so, we introduce the time-varying state transformation
from $e_{s}$ to $\tilde{e}_{s}$ as
\begin{equation}
\tilde{e}_{s}=k_{1}^{-1}\alpha_{s}(\mu)(\tilde{K}\t\Phi(\mu)\otimes I_{n})e_{s},\label{eq:til_e_s-1}
\end{equation}
which coincides with the procedure in Section \ref{sec:A-Cascade-Design}.

Define functions $B(s)=k_{1}^{-1}(\alpha_{x}(s))^{-L_{m}}$, $\delta_{s}(s)=\frac{\mathrm{d}\alpha_{s}(s)}{\mathrm{d}s}s^{2}(\alpha_{s}(s))^{-1}$ and $ \delta_{x}(s)=\frac{\mathrm{d}\alpha_{x}(s)}{\mathrm{d}s}s^{2}\left(\alpha_{x}(s)\right)^{-1}$.
By (\ref{eq:chain_sys}), (\ref{eq:til_s^i-2}), and (\ref{eq:til_e_s^i}),
$\tilde{e}_{s}$-dynamics can be expressed as
\begin{equation}
\dot{\tilde{e}}_{s}
=\alpha_{s}(\mu)(B(\mu)(u+\varphi(x,d)+\pi(x))-\dot{\varpi}^{i}+\delta_{s}(\mu)\tilde{s}),\label{eq:dot_omega-2}
\end{equation}
where $\pi(x)=\alpha_{x}(\mu)^{L_{m}}K\t\dot{r}_{1}-L_{m}\delta_{x}(\mu)x_{m}$
and
\begin{align}
\dot{r}_{1}= & \big[x_{2};\left(\alpha_{x}(\mu)\right)^{-L_{2}}(x_{3}-L_{2}\delta_{x}(\mu)x_{2});\nonumber \\
 & \;\cdots;\left(\alpha_{x}(\mu)\right)^{-L_{m-1}}(x_{m}-L_{m-1}\delta_{x}(\mu)x_{m-1})\big].\label{eq:r_2-1}
\end{align}
Since $\Lambda$ is Hurwitz, there exist positive matrices $P$ and $Q$
such that $P(\Lambda\otimes I_{n})+(\Lambda\t\otimes I_{n})P=-Q$.
Define two constants
\begin{equation}
v_{1}=\lambda_{\min}(Q)\lambda_{\max}^{-1}(P),\quad v_{2}=2m\lambda_{\max}(P)\lambda_{\min}^{-1}(P).\label{eq:c^i_1}
\end{equation}
Then, we propose the following design criteria \textbf{(DC)}
for $\mathcal{K}_{\infty}$ functions $\alpha_{x}(s)$ in (\ref{eq:til_s^i-2})
and $\alpha_{s}(s)$ in (\ref{eq:til_e_s-1}) such that the time-varying
state transformation (\ref{eq:til_e_s-1}) and the $\tilde{e}_{s}$-dynamics
satisfy $\textbf{C}_{\textbf{2}}'$ in Section \ref{subsec:TST}.

$\textbf{DC}_{\textbf{1}}$: $\alpha_{x}(s)$ satisfies
$\frac{\mathrm{d}\alpha_{x}(s)}{\mathrm{d}s}\leq\frac{v_{1}}{2v_{2}}s^{-2}\left(\alpha_{x}(s)\right)^{2}$ and $\alpha_{x}(s)\leq\frac{c^{*}}{v_{1}}\alpha(s)$,
where $v_{1}$, $v_{2}$ are given in (\ref{eq:c^i_1}) and $c^{*}$
is given in (\ref{eq:c_1});

$\textbf{DC}_{\textbf{2}}$: $\alpha_{s}(s)$ is chosen as
 $\alpha_{s}(s)=(\alpha_{x}(s))^{m}\exp\big(\frac{v_{1}}{2}\sint_{0}^{s}\tau^{-2}\alpha_{x}(\tau)\mathrm{d}\tau\big)$.

\blem \label{lem:5-1}Consider the system (\ref{eq:chain_sys}),
$\varpi^{i}$-dynamics in (\ref{eq:dot_varpi_i-1}) and $p^{i}$-dynamics
in (\ref{eq:dot_p_i}) with time-varying state transformation (\ref{eq:til_e_s-1}).
If conditions in Theorem \ref{the:auxi_genera} and two design criteria
$\textbf{DC}_{\textbf{1}}$-$\textbf{DC}_{\textbf{2}}$ hold, then
the bound of $e_{s}$ satisfies
\begin{equation}
\Vert e_{s}\Vert\leq\tilde{\varepsilon}_{s}^{e}(\|\tilde{e}_{s}\|_{\mathcal{T}})\kappa^{-\frac{v_{1}}{4m}}\left(\alpha_{x}(\mu)\right)\label{eq:norm_e_s}
\end{equation}
 for some $\mathcal{K}_{\infty}^{e}$ function $\tilde{\varepsilon}_{s}^{e}$
and $\|\tilde{e}_{s}\|_{\mathcal{T}}=\sup_{t\in\mathcal{T}_{p}}\|\tilde{e}_{s}(t)\|$.
\elem

As provided in the appendix, the proof of Lemma \ref{lem:5-1}
implies that when $\tilde{e}_{s}$ is bounded for $t\in\mathcal{T}_{p}$,
the prescribed-time convergence of $e_{s}$ is obtained. Therefore,
it suffices to design the controller $u$ in (\ref{eq:dot_omega-2})
such that the closed-loop system for $\tilde{e}_{s}$ admits a prescribed-time
ISS Lyapunov function as in \textbf{$\textbf{C}_{\textbf{2}}'$} and $\tilde{e}_{s}$
is bounded for $t\in\mathcal{T}_{p}$. Then, we design the controller
$u$ as
\begin{align}
u= & -(v+\psi^{2}(x)+1)\mbox{sign}(k_{1})\tilde{e}_{s}-\pi(x)\nonumber \\
 & -(B(\mu))^{-1}\delta_{s}(\mu)\tilde{s}\label{eq:u^i}
\end{align}
 with $v>0$, and function $\psi$ defined in (\ref{eq:psi^i}).

For simplicity, we define
\begin{equation}
\alpha_{b}(s)=\alpha_{s}(s)|B(s)|,\quad\tilde{\alpha}_{b}(s)=\alpha_{s}(s)|B(s)|^{-1},\label{eq:alpha_b}
\end{equation}
 where $\alpha_{s}(s)$ is
introduced in (\ref{eq:til_e_s^i}). We present the following proposition, whose proof is given in appendix.
\bproposition \label{prop:5-1}
For $s\in \mathbb R_{\geq 0}$, $\alpha_b(s)$ and $\tilde \alpha_b(s)$ are $\mathcal{K}_{\infty}$ functions.
\eproposition

\blem \label{lem:5-2} Consider the system (\ref{eq:chain_sys})
with the controller (\ref{eq:u^i}), $\varpi^{i}$-dynamics in (\ref{eq:dot_varpi_i-1})
and $p^{i}$-dynamics in (\ref{eq:dot_p_i}) with time-varying state
transformation (\ref{eq:til_e_s-1}). If conditions in Theorem \ref{the:auxi_genera}
and two design criteria $\textbf{DC}_{\textbf{1}}$-$\textbf{DC}_{\textbf{2}}$
hold, then $\tilde{e}_{s}$-dynamics satisfies condition $\textbf{C}_{\textbf{2}}'$.
Moreover, it admits the prescribed-time ISS Lyapunov function in \textbf{$\textbf{C}_{\textbf{2}}'$} (omitting superscript $i$) with
\begin{equation}
\begin{gathered}\underline{\alpha}_{\tilde{s}}(s)=\bar{\alpha}_{\tilde{s}}(s)=s/2,\quad\tilde{\alpha}_{\tilde{s}}(s)=2v\alpha_{b}(s),\\
\sigma_{\tilde{d}}(s)=h^{2}(s),\quad\tilde{\sigma}_{\tilde{d}}(s)=\alpha_{b}(s)/4,\\
\sigma_{\tilde{r}}(s)=s^{2},\quad\tilde{\sigma}_{\tilde{r}}(s)=\tilde{\alpha}_{b}(s)(\gamma_{\zeta}(s))^{2}/4,
\end{gathered}
\label{eq:con_1}
\end{equation}
where $\gamma_{\zeta}(s)=\max\{2\lambda_{N}+\varrho_{c},1\}\alpha(s)$.
And the controller $u$ satisfies $\text{\textbf{C}}_{\textbf{3}}'$
with
\begin{equation}
\tilde{\gamma}_{s}(s)=\varepsilon_{1}'+\varepsilon_{2}'s,\quad\tilde{\gamma}_{u}(s)=\varepsilon_{3}'(\alpha_{x}(s))^{m}\label{eq:til_gamma_s-1}
\end{equation}
for some finite constants $\varepsilon_{1}'$, $\varepsilon_{2}'$
and $\varepsilon_{3}'$. \elem

Applying Theorems \ref{thm:2}, \ref{the:auxi_genera} and Lemmas \ref{lem:5-1},
\ref{lem:5-2}, we obtain the following results.

\bthm\label{thm:5-1} Consider the system composed of (\ref{eq:dot_varpi_i-1}),
(\ref{eq:dot_p_i}), (\ref{eq:chain_sys}) and (\ref{eq:u^i}). If
conditions in Theorem \ref{the:auxi_genera} and two design criteria
$\textbf{DC}_{\textbf{1}}$-$\textbf{DC}_{\textbf{2}}$ hold, the
DPTCO problem for the chain integrator MASs (\ref{eq:chain_sys})
is solved. \ethm

\section{Adaptive DPTCO for Strict-Feedback MASs\label{sec:Adaptive-DPTCO-for} }
In this section, to further examine the generality of proposed DPTCO framework proposed in Section \ref{sec:A-Cascade-Design}, 
we consider the adaptive DPTCO problem for a class of nonlinear strict-feedback MASs with parameter uncertainty, as follows,   
\begin{equation}
\begin{aligned}\text{\ensuremath{\dot{x}_{1}}} & =x_{2},\\
\dot{x}_{q} & =x_{q+1}+\theta\varphi_{q}(x_{q}),\;q=2,\cdots,m-1,\\
\dot{x}_{m} & =u+\theta\varphi_{m}(x_{m}),\quad y  =x_{1},\\
\end{aligned}
\label{eq:strict-feedback system}
\end{equation}
 where $x=[x_{1}\t,\cdots,x_{m}\t]\t\in\mathbb{R}^{mn}$ is the system
state with $x_{q}\in\mathcal{\mathbb{R}}^{n}$, $y\in\mathbb{R}^{n}$
is output and $u\in\mathbb{R}^{n}$ is control input. $\theta\in\mathbb{R}$
is an unknown constant and $\varphi_{q}(x_{q}):\mathbb{R}^{n}\to\mathbb{R}^{n}$
is a known function with $\varphi_{q}(0)=0$ for $q=2,\cdots,m$.
\textit{For simplicity, we omit the superscript $i$ when no confusion
is raised. }

\bass \label{ass:varphi-1} \cite{song2016adaptive44}
For $q=2,\cdots,m$, $\varphi_{q}$ is first-order differentiable
and locally Lipschitz function.\eass

\brem \label{rem:6-1}Under Assumption \ref{ass:varphi-1}, and given that  $\varphi_{q}(0)=0$, by the mean value theorem, there exists a 
continuous matrix-valued function $\psi_{q}(x_{q}):\mathbb{R}^{n}\to\mathbb{R}^{n\times n}$
such that
\begin{equation}
\varphi_{q}(x_{q})=\psi_{q}(x_{q})x_{q},\label{eq:Hurwitz}
\end{equation}
 where $\psi_{q}(x_{q})$ and its first derivative with respect to
$t$ are continuous and bounded. Without losing generality, we assume
$\Vert\psi_{q}(x_{q})\Vert\leq\bar{\psi}_{q}$ and $\Vert\frac{\mathrm{d}\psi_{q}(x_{q}(t))}{\mathrm{d}t}\Vert\leq\tilde{\psi}_{q}$
hold for $x_{q}\in\mathbb{R}^{n}$, where $\bar{\psi}_{q}$ and $\tilde{\psi}_{q}$
are some positive finite constants. \erem

Following the procedure in Section \ref{sec:A-Cascade-Design},
we define the error states according to (\ref{eq:e_y_varpi}) as
\begin{equation}
e_{s}=\big[(x_{1}-\varpi^{i})\t,x_{2}\t,\cdots,x_{m}\t,\varsigma\t\big]\t,\label{eq:e_s^i}
\end{equation}
where $\varpi^{i}$ is given in (\ref{eq:dot_varpi_i-1}) and
$
\varsigma=[\hat{\theta},\xi_{f}\t]\t
$
is the controller state where $\hat{\theta}\in\mathbb{R}$ is the
estimator of unknown parameter $\theta$ and $\xi_{f}=[\xi_{2f}\t,\cdots,\xi_{mf}\t]\t\in\mathbb{R}^{n(m-1)}$
is the dynamic filter variable to be designed.
To facilitate the stability analysis and simplify the derivation,
we introduce the coordinate transformation as
\begin{equation}
\begin{aligned} & \tilde{x}_{1}=x_{1}-\varpi^{i},\quad\xi_{1}=-c_{1}\alpha_{\xi}(\mu)\tilde{x}_{1},\\
 & \tilde{x}_{q}=x_{q}-\xi_{qf},\quad\tilde{\xi}_{q}=\xi_{qf}-\xi_{q-1},\\
 & \xi_{q}=-c_{q}\alpha_{\xi}(\mu)\tilde{x}_{q}-\hat{\theta}\varphi_{q}(x_{q})-\upsilon_{q}\alpha_{\xi}(\mu)\tilde{\xi}_{q},\\
 & \qquad\qquad q=2,\cdots m,
\end{aligned}
\label{eq:til_x}
\end{equation}
where $c_{q}$ for $q=1,\cdots,m$ is to be determined, $\alpha_{\xi}\in\mathcal{K}_{\infty}$
to be designed, $\xi=[\xi_{1}\t,\cdots\xi_{m}\t]\t$ is the virtual
controller and $\xi_{qf}$ and $\hat{\theta}$-dynamics are designed
as
\begin{align}
\dot{\xi}_{qf} & =\upsilon_{q}\alpha_{\xi}(\mu)(-\xi_{qf}+\xi_{q-1}),\;q=2,\cdots,m,\label{eq:dot_xi_qf}\\
\dot{\hat{\theta}} & =\tau-\sigma\alpha_{\xi}(\mu)\hat{\theta},\label{eq:adaptive_law}
\end{align}
with $\upsilon_{q}$ for $q=2,\cdots,m$ and $\sigma>0$ to be determined,
and
\begin{equation}
\tau=\ssum_{q=2}^{m}\tau_{q},\quad\tau_{q}=\left(\alpha_{\xi}(\mu)\right)^{2L_{q}}\tilde{x}_{q}\t\varphi_{q}(x_{q}).\label{eq:tau_q^i}
\end{equation}

We further introduce the time-varying state transformation for (\ref{eq:e_s^i})
as
$
\tilde{e}_{s}=\big[\omega\t,\eta\t,\tilde{\theta}\big]\t$,
 where $\omega=\left[\omega_{1}\t,\cdots,\omega_{m}\t\right]\t$,
$\eta=\left[\eta_{2}\t,\cdots,\eta_{m}\t\right]\t$ and $\tilde{\theta}=\theta-\hat{\theta}$
with
\begin{equation}\label{eq:omega eta}
\begin{gathered}
  \omega_{q}=\left(\alpha_{\xi}(\mu)\right)^{L_{q}}\tilde{x}_{q},\;q=1,\cdots,m, \\
\eta_{q}=\left(\alpha_{\xi}(\mu)\right)^{L_{q}}\tilde{\xi}_{q},\;q=2,\cdots,m,
\end{gathered}
\end{equation}
where $L_{q}=m+l+1-q$ with $l>0$,
$q=1,\cdots m$. By doing so, we in fact introduce the time-varying
state transformation from $e_{s}$ to $\tilde{e}_{s}$ as
\begin{equation}
\begin{gathered}\omega  =(\Phi_{1}(\mu)\otimes I_{n})\Lambda_{1}e_{s},\quad \tilde{\theta}  =\theta-\Lambda_{4}e_{s},\\
\eta  =(\Phi_{2}(\mu)\otimes I_{n})(\Lambda_{2}e_{s}-\Lambda_{3}\xi(e_{s}))
\end{gathered}
\label{eq:omega_eta}
\end{equation}
with
$
\Phi_{1}(\mu)=\mbox{diag}\{(\alpha_{\xi}(\mu))^{L_{1}},\cdots,(\alpha_{\xi}(\mu))^{L_{m}}\}$,
$ 
\Phi_{2}(\mu)=\mbox{diag}\{(\alpha_{\xi}(\mu))^{L_{2}},\cdots,(\alpha_{\xi}(\mu))^{L_{m}}\}$, 
and 
 $\Lambda_{1} = \big[I_{nm},0_{nm\times1},
(0_{n\times n(m-1)}\t, I_{n(m-1)})\t
\big]$, 
$\Lambda_{2}=[0_{n(m-1)\times(nm+1)},I_{n(m-1)}]$, $\Lambda_{3}=[I_{n(m-1)},0_{n(m-1)\times n}]$ and 
$\Lambda_{4}=[0_{1\times nm},1,0_{1\times n(m-1)}]$. 
As a result, the $\tilde{e}_{s}$-dynamics can be expressed as $\dot{\tilde{e}}_{s}=\tilde{h}_{s}(\tilde{e}_{s},e_{r},u,\theta,\mu).$
We propose the design criterion for $\mathcal{K}_{\infty}$ functions
$\alpha(s)$ and $\alpha_{\xi}(s)$.

$\textbf{DC}_{\xi}$: $\alpha_{\xi}(s)$ satisfies $\frac{\mathrm{d}\alpha_{\xi}(s)}{\mathrm{d}s}\leq s^{-2}(\alpha_{\xi}(s))^{2}$ and $\alpha_{\xi}(s)\leq\frac{c^{*}\alpha(s)}{2L_{2}}$ 
for $s\in\mathbb{R}_{\geq0}$, where $c^{*}$ is denoted in
(\ref{eq:c_1}).

\blem \label{lem:6-1} Consider the system (\ref{eq:strict-feedback system}),
$\varpi^{i}$-dynamics in (\ref{eq:dot_varpi_i-1}) and $p^{i}$-dynamics
in (\ref{eq:dot_p_i}) with time-varying state transformation (\ref{eq:omega_eta}).
If conditions in Theorem \ref{the:auxi_genera} and the design criterion
$\textbf{DC}_{\xi}$ hold, then the bound of $e_{s}$
satisfies
\begin{equation}
\Vert e_{s}\Vert\leq\varepsilon_{s}^{e}(\|\tilde{e}_{s}\|_{\mathcal{T}})(\alpha_{\xi}(\mu))^{-1}\label{eq:norm_e_s-2}
\end{equation}
 for some $\mathcal{K}_{\infty}^{e}$ function $\varepsilon_{s}^{e}$
and $\|\tilde{e}_{s}\|_{\mathcal{T}}=\sup_{t\in\mathcal{T}_{p}}\|\tilde{e}_{s}(t)\|$.
\elem

The proof of Lemma \ref{lem:6-1} is given in the appendix. It implies
that when $\tilde{e}_{s}$ is bounded for $t\in\mathcal{T}_{p}$,
the prescribed-time convergence of $e_{s}$ is achieved. Then,
the controller $u$ is designed as 
\begin{equation}
u=\xi_{m}\label{eq:u^i-1}
\end{equation}
where $\xi_{m}$ is designed in (\ref{eq:til_x}).

\bthm \label{thm:6-2} Consider the system (\ref{eq:strict-feedback system})
with the controller (\ref{eq:u^i-1}), $\varpi^{i}$-dynamics in (\ref{eq:dot_varpi_i-1})
and $p^{i}$-dynamics in (\ref{eq:dot_p_i}) with time-varying state
transformation (\ref{eq:omega_eta}) under Assumption \ref{ass:varphi-1}.
Suppose conditions in Theorem \ref{the:auxi_genera} and the
design criterion $\textbf{DC}_{\xi}$ hold. Then,
there always exists a set of parameters $c_{q}$ for $q=1,\cdots,m$,
$\upsilon_{q}$ for $q=2,\cdots,m$ and $h$ such that $\Omega(h)$
is an invariant set where $\Omega(h)=\{\tilde{e}_{s}\in\mathbb{R}^{2mn+1-n}|\|\tilde{e}_{s}\|^{2}\leq h^{2}\}$
and the DPTCO problem for strict-feedback MASs (\ref{eq:strict-feedback system})
is solved. \ethm

\section{Simulation Results \label{sec:Simulation-Results}}

In this section, we show two numerical examples to illustrate the
theoretical results. The graph for the two simulations is given by $1\leftrightarrow 2\leftrightarrow3\leftrightarrow4\leftrightarrow5 \leftrightarrow6 \leftrightarrow1$.
\begin{example} (Robust DPTCO for Euler-Lagrange MASs)
Consider the Euler-Lagrange MASs as $\dot{x}_{1}^{i}=  x_{2}^{i}$, $\dot{x}_{2}^{i}=  M(x_{1}^{i})^{i,-1}(u^{i}-C^{i}(x_{1}^{i},x_{2}^{i})x_{2}^{i}-G^{i}(x_{q}^{i}))$, $y^{i}=  x_{1}^{i},i=1,\cdots,6$
where $x_{1}^{i},x_{2}^{i}\in\mathbb{R}^{2}$ with $x_{1}^{i}=[x_{11}^{i},x_{12}^{i}]\t$,
$x_{2}^{i}=[x_{21}^{i},x_{22}^{i}]\t$, $M^{i}(x_{1}^{i})=[\theta_{1}^{i}+\theta_{2}^{i}+2\theta_{3}^{i}\cos(x_{12}^{i}),\theta_{2}^{i}+\theta_{3}^{i}\cos(x_{12}^{i});\theta_{2}^{i}+\theta_{3}^{i}\cos(x_{12}^{i}),\theta_{4}^{i}  ]$, $C^{i}(x_{1}^{i},x_{2}^{i})=[-\theta_{3}^{i}\sin(x_{12}^{i})x_{21}^{i}, -2\theta_{3}^{i}\sin(x_{12}^{i})x_{21}^{i}; 0, \theta_{3}^{i}\sin(x_{12}^{i})x_{22}^{i} ]$ and $G^{i}(x_{q}^{i})= [\theta_{5}^{i}g\cos(x_{11}^{i})+\theta_{6}^{i}g\cos(x_{11}^{i}+x_{12}^{i});\theta_{6}^{i}g\cos(x_{11}^{i}+x_{12}^{i})]$. The unknown parameters are 
$\theta_{1}^{i}=7$, $\theta_{2}^{i}=0.96$, $\theta_{3}^{i}=1.2$,
$\theta_{4}^{i}=5.96$, $\theta_{5}^{i}=2$, $\theta_{6}^{i}=1.2$
 for $i=1,\cdots,6$, and $g=9.8$. Note that
the system is in the form of the chain-integrator systems in (\ref{eq:chain_sys})
and satisfies (\ref{eq:psi^i}) due to the structural property of
Euler-Lagrange systems.

The six robots are located in a thermal radiation field, and the relationship
between the intensity of thermal radiation $P$, temperature $T_{em}$
and distance $d$ can be roughly expressed as
$
P\propto\frac{T_{em}^{4}}{\|d-d^{*}\|^{2}}
$, 
where $d^{*}$ denotes the two-dimensional coordinates of the heat
source. Suppose each robot is capable of measuring the gradient information
of the heat source with respect to distance. The objective is to design
controller $u^{i}$ such that the six robots approach the heat source
in a formation, and reduce the total displacement of the six robots
from their original location. Thus, the global objective function
is designed as
$\min\ssum_{i=1}^{6}\iota_{1}^{i}\|y^{i}-d^{*}\|^{2}+\ssum_{i=1}^{6}\iota_{2}^{i}\|y^{i}-y^{i}(t_{0})\|^{2},\; \mbox{s.t. }\;y^{i}-y^{j}=\omega^{i}-\omega^{j}$
where $\omega^{1}=[1,0]\t$, $\omega^{2}=[1/2,\sqrt{3}/2]\t$, $\omega^{3}=[-1/2,\sqrt{3}/2]\t$,
$\omega^{4}=[-1,0]\t$, $\omega^{5}=[-1/2,-\sqrt{3}/2]\t$, $\omega^{6}=[1/2,-\sqrt{3}/2]\t$
represent the formation shape, and $\iota_{1}^{i}$ and $\iota_{2}^{i}$
are objective weights.
By defining $\bar{y}^{i}=y^{i}-\omega^{i}$, the optimization problem
 is transformed into $\min\ssum_{i=1}^{6}\iota_{1}^{i}\|\bar{y}^{i}+\omega^{i}-d^{*}\|^{2}+\ssum_{i=1}^{6}\iota_{2}^{i}\|\bar{y}^{i}+\omega^{i}-y^{i}(t_{0})\|^{2},\;\mbox{s.t.}\;\bar{y}^{i}=\bar{y}^{j} $, 
which is consistent with (\ref{eq:opti_problem}). For the optimization
problem, we design $\zeta^{i}$-dynamics as
in the form of (\ref{eq:dot_varpi_i-1}), (\ref{eq:dot_p_i}) such
that $\varpi^{i}$ converges to the optimum within prescribed
time. Then, the reference trajectory for each robot dynamics is changed
as
$
\varpi^{i,'}=\varpi^{i}+\omega^{i}
$. 
Replacing $\varpi^{i}$ in Section \ref{sec:Optimal-Trajectory-Generator}
with $\varpi^{i,'}$, we can design the controller following the procedures
in Section \ref{sec:DPTCO-for-Chain} to solve the optimization problem. Let the initial condition  be $x_{1}^{1}(t_{0})=[1,1]\t$,
$x_{1}^{2}(t_{0})=[2,2]\t$, $x_{1}^{3}(t_{0})=[3,3]\t$, $x_{1}^{4}(t_{0})=[-2,-3]\t$,
$x_{1}^{5}(t_{0})=[-2,-1]\t$, $x_{1}^{6}(t_{0})=[0,0]\t$, $x_{2}^{i}(t_{0})=[1,1]\t$,
$\varpi^{i}(t_{0})=[1,1]\t$, $p^{i}(t_{0})=[0,0]\t$ for $i=1,\cdots,6$.
The initial time $t_{0}$ is set as $t_{0}=0$, and the prescribed-time
scale $T=1s$. The parameters and gain functions are chosen as $c=6$,
$l=2$, $k_{1}=2$, $\alpha(\mu)=10\mu$, $\alpha_{x}(\mu)=\mu$,
$\alpha_{s}(\mu)=\exp(\mu)$. The weight coefficients $\iota_{1}^{i}$
and $\iota_{2}^{i}$ for objective function
are chosen as $\iota_{1}^{i}=0.5$ and $\iota_{2}^{i}=0.1$ for $i=1,\cdots,6$.
The coordinate of heat source is set as $d^{*}=[0,0]\t$.

\begin{figure}[htbp]
	\centering
	\includegraphics[width=0.8\linewidth]{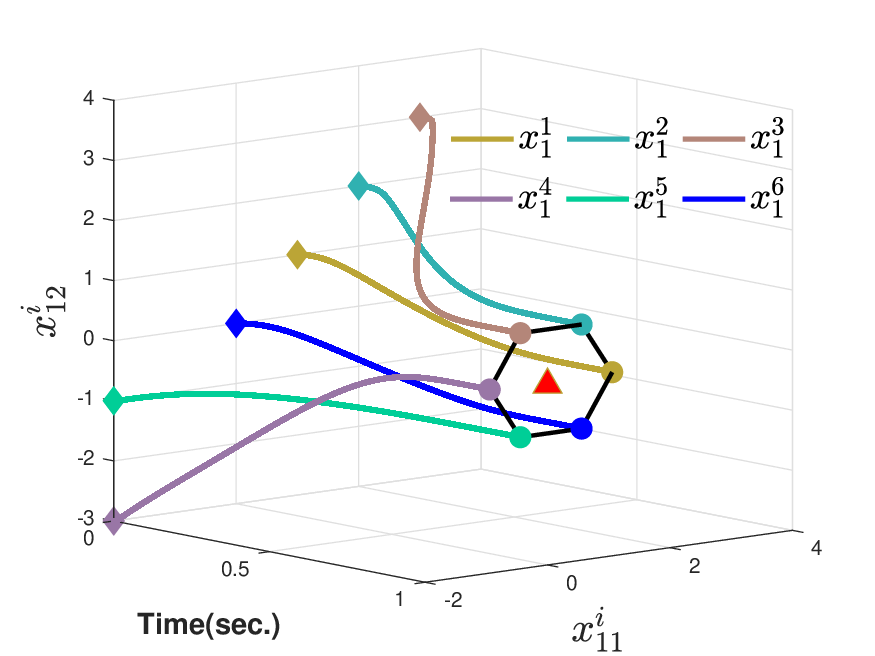}
	\caption{ Trajectories of positions $x_{1}^{i}$ of the six
robots for $0\protect\leq t<T$, where  $\blacklozenge$ and $\bullet$
denote the initial and final positions, respectively,  and $\blacktriangle$ denotes the heat source. }
\label{fig:error}
\end{figure}

\begin{figure}[htbp]
	\centering
	\includegraphics[width=0.8\linewidth]{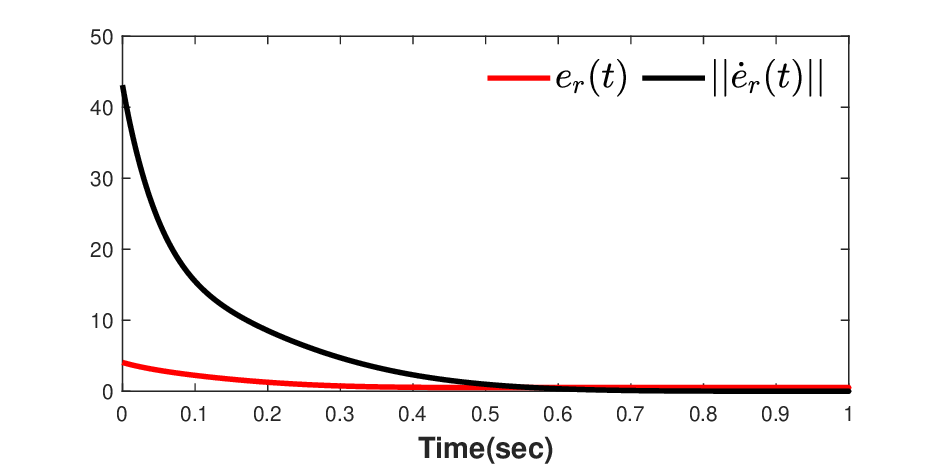}
	\caption{ The trajectories of $e_{r}(t)$ and $\dot{e}_{r}(t)$. }
\label{fig:dot_varpi_p}
\end{figure}

\begin{figure}[htbp]
	\centering
	\includegraphics[width=0.8\linewidth]{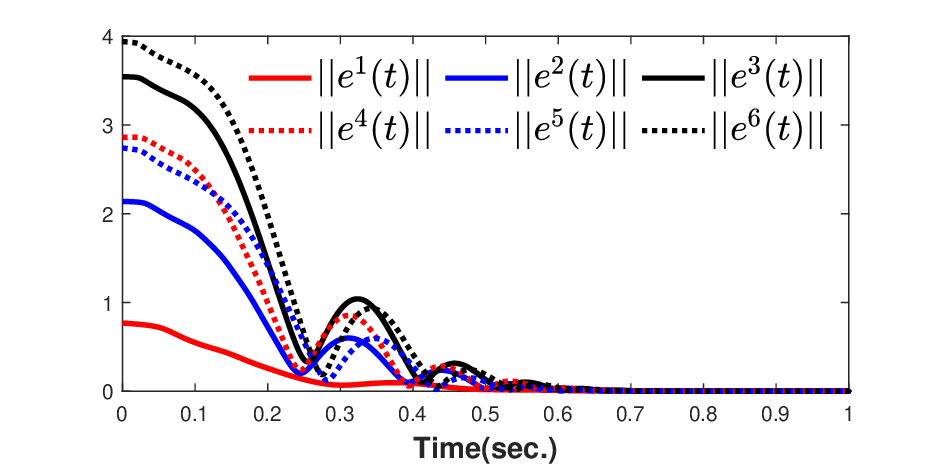}
	\caption{ The trajectories of tracking errors between
	each agent's output and optimum. }
\label{fig:error_SFS}
\end{figure}
\begin{figure}[htbp]
	\centering
	\includegraphics[width=0.8\linewidth]{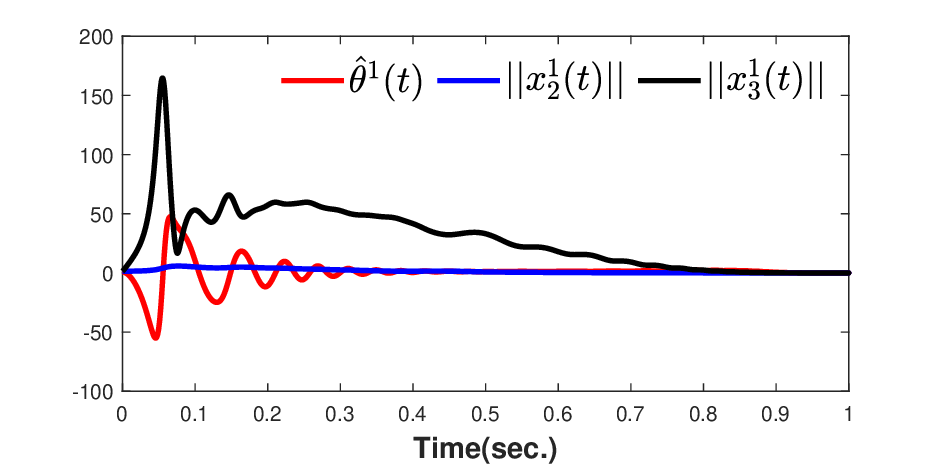}
	\caption{ The trajectories of $\hat{\theta}^{1}(t)$,
		$\Vert x_{2}^{1}(t)\Vert$, and $\Vert x_{3}^{1}(t)\Vert$.}
\label{fig:hat_theta}
\end{figure}

The simulation results are shown in Figures.
 \ref{fig:error} and \ref{fig:dot_varpi_p}. In Figure. \ref{fig:dot_varpi_p}, $e_{r}(t)$
and $\dot{e}_{r}(t)$ converge to zero within $T$, and thus the validity
of the optimal trajectory generator designed in Section \ref{sec:Optimal-Trajectory-Generator}
is verified. In Figure. \ref{fig:error}, the six robots approach
the heat source in formation within the prescribed time. \end{example}

\begin{example} (Adaptive DPTCO for strict-feedback MASs) Consider
the strict-feedback MASs in the presence of parameter uncertainties
as $\dot{x}_{1}^{i}=x_{2}^{i}$, $\dot{x}_{2}^{i}=x_{3}^{i}+\theta^{i}x_{2}^{i}$, $\dot{x}_{3}^{i}=u^{i}+\theta^{i}x_{3}^{i}$, $y^{i}=x_{1}^{i},i=1,\cdots,6$, 
where $\theta^{i}\in\mathbb{R}$, $x_{1}^{i},x_{2}^{i},x_{3}^{i}\in\mathbb{R}^{2}$.
$\theta=[\theta^{1},\cdots,\theta^{6}]=[1,2,-1,3,-2,-3]$.  The
local objective function of each agent is $f^{i}(z)=  \exp\left((z-z_{d2}^{i})\t P^{i}(z-z_{d2}^{i})\right)+(z-z_{d1}^{i})\t Q^{i}(z-z_{d1}^{i})+\delta^{i}$, 
where $z_{d1}=[z_{d1}^{1},\cdots,z_{d1}^{6}]=1_{2\times6}$, $z_{d2}=[z_{d2}^{1},\cdots,z_{d2}^{6}]=0.5_{2\times6}$,
$\delta=[\delta^{1},\cdots,\delta^{6}]=1_{6}$, 
$P^i$ and $Q^i$ are positive definite matrices. Using Global Optimization
Toolbox in MATLAB, the optimal agreement $z^{*}$ is
$
z^{*}=[0.7263,0.7183]\t
$, 
which is used for verification only. The parameters are chosen as
$c_{1}^{i}=10$, $c_{2}^{i}=10$, $c_{3}^{i}=10$, $\upsilon_{2}^{i}=15$,
$\upsilon_{3}^{i}=20$, $\sigma^{i}=10,i=1,\cdots,6$, $l=1$. $\alpha(\mu)=10\mu^{3/2}$,
$\alpha_{\xi}(\mu)=\mu^{3/2}$. The initial values are $x_{1}^{1}(t_{0})=[1,0]\t$,
$x_{1}^{2}(t_{0})=[2,-1]\t$, $x_{1}^{3}(t_{0})=[3,-2]\t$, $x_{1}^{4}(t_{0})=[-1,3]\t$,
$x_{1}^{5}=[-2,1]\t$, $x_{1}^{6}(t_{0})=[-3,2]\t$, $x_{2}^{i}(t_{0})=[1,1]\t$,
$x_{3}^{i}(t_{0})=[1,1]\t$, $\hat{\theta}^{i}(t_{0})=1,i=1,\cdots,6$,
and $\varpi(t_{0})$, $p(t_{0})$ are the same as that in Example
1.
The simulation results are shown in Figure. \ref{fig:error_SFS} and
\ref{fig:hat_theta}. In Figure. \ref{fig:error_SFS}, the tracking
error between each agent's output and optimum is bounded
and achieves prescribed-time convergence towards zero. For simplicity,
we only provide the trajectories of $\hat{\theta}^{1}(t)$, $\Vert x_{2}^{1}(t)\Vert$,
$\Vert x_{3}^{1}(t)\Vert$ in Figure. \ref{fig:hat_theta}. These
trajectories show that we achieve prescribed-time convergence towards
zero for $\hat{\theta}^{1}(t)$, $\Vert x_{2}^{1}(t)\Vert$ and $\Vert x_{3}^{1}(t)\Vert$.
\end{example}

\section{Conclusion \label{sec:Conclusion}}
In this paper, we propose a novel DPTCO algorithm for a class of high-order
nonlinear MASs. A DPTCO framework is first constructed embedding the
cascade design such that the DPTCO problem is divided into optimum seeking for thewhole system and reference trajectory tracking
problem for each agent. The DPTCO framework is then utilized to solve
DPTCO problem for chain integrator MASs and strict-feedback MASs.
The prescribed-time convergence lies in the time-varying gains which
increase to infinity as time approaches the prescribed time. When
solving the tracking problem for the two specific MASs, high-order
derivative of reference trajectory is not required. It would be very
interesting to further consider the DPTCO where the local objective functions
subject to bound, equality, and inequality constraints.

\section*{References}
\bibliography{ref2.bib}
\bibliographystyle{ieeetr}

\appendix{}

\prooflater{Theorem \ref{the:auxi_genera}}
 Let us introduce the state transformations as
\begin{equation}
\varphi:=\left[\begin{array}{c}
\bar{\varphi}\\
\tilde{\varphi}
\end{array}\right]=\left[\begin{array}{c}
\bar{r}\t\\
\bar{R}\t
\end{array}\right]e_{\varpi},\;\phi:=\left[\begin{array}{c}
\bar{\phi}\\
\tilde{\phi}
\end{array}\right]=\left[\begin{array}{c}
\bar{r}\t\\
\bar{R}\t
\end{array}\right]e_{p},\label{eq:st}
\end{equation}
where $\bar{\varphi},\bar{\phi}\in\mathbb{R}^{m}$ and $\tilde{\varphi},\tilde{\phi}\in\mathbb{R}^{(N-1)m}$.
As a result, the system composed of (\ref{eq:dot_e_varpi}) and (\ref{eq:dot_e_p})
can be rewritten as 
\begin{equation}
\begin{aligned}\dot{\bar{\varphi}} & =-\alpha(\mu)\bar{r}\t\nabla\tilde{F}(e_{\varpi}),\\
\dot{\tilde{\varphi}} & =-\alpha(\mu)(\bar{\mathcal{L}}_{R}\tilde{\varphi}+\bar{R}\t\nabla\tilde{F}(e_{\varpi})+\tilde{\phi}),\\
\dot{\bar{\phi}} & =0,\quad\dot{\tilde{\phi}}=\alpha(\mu)\bar{\mathcal{L}}_{R}\tilde{\varphi},
\end{aligned}
\label{eq:dot_phi_2}
\end{equation}
 where we used $\bar{r}\t\mathcal{L}=0$ and the fact $\bar r \t e_p = \bar r\t \left( p + \nabla F(1_N \otimes z^*)\right) 
  = \ssum_{i=1}^N p^i + \ssum_{i=1}^N f^i (z^*) =0$. 
Let the Lyapunov function candidate be 
\begin{equation}
V(\varphi,\phi)=\frac{c_{1}}{2}(\varphi\t\varphi+\phi\t\tilde{\mathcal{L}}_{R}^{-1}\phi)+\frac{1}{2}(\varphi+\phi)\t(\varphi+\phi)\label{eq:V}
\end{equation}
where $\tilde{\mathcal{L}}_{R}=\mbox{diag}\{I_{m},\bar{\mathcal{L}}_{R}\}$.
Then the time derivative of $V(\varphi,\phi)$ along (\ref{eq:dot_phi_2})
is $\dot{V}(\varphi,\phi)=c_{1}\alpha(\mu)\big(-\tilde{\varphi}\t\bar{\mathcal{L}}_{R}\tilde{\varphi}-\bar{\varphi}\t\bar{r}\t\nabla\tilde{F}(e_{\varpi})-\tilde{\varphi}\t\bar{R}\t\nabla\tilde{F}(e_{\varpi})\big)+\alpha(\mu)\big(-\bar{\varphi}\t\bar{r}\t\nabla\tilde{F}(e_{\varpi})-\tilde{\varphi}\t\bar{R}\t\nabla\tilde{F}(e_{\varpi})-\tilde{\varphi}\t\tilde{\phi}-\tilde{\phi}\t\tilde{\phi}-\bar{\phi}\t\bar{r}\t\nabla\tilde{F}(e_{\varpi})-\tilde{\phi}\t\bar{R}\t\nabla\tilde{F}(e_{\varpi})\big)$,
where we used $\phi\t\dot{\phi}=\tilde{\phi}\t\dot{\tilde{\phi}}+\bar{\phi}\t\dot{\bar{\phi}}$
and $\bar{\phi}\t\dot{\bar{\phi}}=0$. Due to Assumption \ref{ass:cost_func}
and (\ref{eq:st}), one has $-\bar{\varphi}\t\bar{r}\t\nabla\tilde{F}(e_{\varpi})-\tilde{\varphi}\t\bar{R}\t\nabla\tilde{F}(e_{\varpi})\leq-\rho_{c}\Vert e_{\varpi}\Vert^{2}=-\rho_{c}\Vert\varphi\Vert^{2}$. 
A few facts are $\bar{\phi}\t\bar{\phi}=0$, $-(\bar{\phi}\t\bar{r}\t+\tilde{\phi}\t\bar{R}\t)\nabla\tilde{F}(e_{\varpi})\leq\|\phi\|^{2}/4+\varrho_{c}^{2}\|\varphi\|^{2}$
and $-\tilde{\varphi}\t\tilde{\phi}\leq\|\tilde{\varphi}\|^{2}+\|\phi\|^{2}/4$.
Substituting the above inequalities into $\dot{V}(\varphi,\phi)$
yields
\begin{align}
\dot{V}(\varphi,\phi) & \leq-\alpha(\mu)(c_{1}\rho_{c}+\rho_{c}-\varrho_{c}^{2})\|\varphi\|^{2}\nonumber \\
 & \quad-\alpha(\mu)(c_{1}\lambda_{2}-1)\|\tilde{\varphi}\|^{2}-\alpha(\mu)\|\phi\|^{2}/2.\label{eq:phi_1}
\end{align}
By (\ref{eq:c_1}), we have $c_{1}\lambda_{2}\geq1$
and $c_{1}\rho_{c}+\rho_{c}-\varrho_{c}^{2}\geq1/2$. Then, by (\ref{eq:phi_1}),
one has
$
\dot{V}(\varphi,\phi)\leq-\alpha(\mu)\Vert\mbox{col}[\varphi,\phi]\Vert^{2}/2
$. 
Note that $V(\varphi,\phi)$ can be written as a function of $e_{r}$, 
i.e., $V(e_{r})=V(\varphi,\phi)=c_{1}(\varphi\t\varphi+\phi\t\tilde{\mathcal{L}}_{R}^{-1}\phi)/2+(\varphi+\phi)\t(\varphi+\phi)/2=c_{1}\big(\|e_{\varpi}\|^{2}+e_{p}\t[r,R]\tilde{\mathcal{L}}_{R}^{-1}[r,R]\t e_{p}\big)/2+\|e_{\varpi}+e_{p}\|^{2}/2$.
As a result,
\begin{equation}
\begin{gathered}c_{2}\Vert e_{r}\Vert^{2}\leq V(e_{r})\leq c_{3}\Vert e_{r}\Vert^{2},\\
\dot{V}(e_{r})\leq-\alpha(\mu)\Vert e_{r}\Vert^{2}/2\leq-2c^{*}\alpha(\mu)V(e_{r}),
\end{gathered}
\label{eq:con3-2}
\end{equation}
where $c^{*}$ is given in (\ref{eq:c_1}). Thus, the first part of
condition \textbf{$\textbf{C}_{\textbf{1}}$} is satisfied. From (\ref{eq:dot_e_varpi})-(\ref{eq:dot_e_p}),
one has
\begin{equation}
\|\dot{e}_{r}\|=\Vert\tilde{h}_{\zeta}\left(e_{r},\mu\right)\Vert\leq\gamma_{\zeta}(\mu)\Vert e_{r}\Vert,\label{eq:til_h_zeta}
\end{equation}
where $\gamma_{\zeta}(\mu)=\max\{2\lambda_{N}+\varrho_{c},1\}\alpha(\mu)$.
Invoking comparison lemma for (\ref{eq:con3-2}) leads to (\ref{eq:e_r_bound})
with $\gamma_{r}(\|e_{r}(t_{0})\|)=\sqrt{c_{3}/c_{2}}\Vert e_{r}(t_{0})\|$.
By (\ref{eq:par_alpha}), one has 
\begin{align}
&\alpha(\mu)\kappa^{-c^{*}}(\alpha(\mu))\notag \\
& = \alpha(\mu(t_0)) \exp \left(\sint_{t_0}^t \left.\frac{\mathrm d\alpha(s)}{\mathrm ds}\right|_{s= \mu(\tau )} \frac{\dot \mu(\tau)}{\alpha(\mu(\tau))}\mathrm d\tau\right)\kappa^{-c^{*}}(\alpha(\mu))\notag \\
&\leq \alpha(\mu(t_{0}))\kappa^{-\frac{c^{*}}{2}}\left(\alpha(\mu)\right)\label{eq:alpha_kappa}
\end{align}
where we used the facts that $\dot \mu(t)=\mu^2 (t)$ and $\frac{\mathrm d\alpha(s)}{\mathrm ds}s^2\leq c^* (\alpha(s))^2/2$. 
Then substituting (\ref{eq:e_r_bound}) into (\ref{eq:til_h_zeta})
and utilizing (\ref{eq:alpha_kappa}) leads to (\ref{eq:er_bound-1})
with $\gamma_{e}(\|e_{r}(t_{0})\|)=\max\{2\lambda_{N}+\varrho_{c},1\}\alpha(\mu(t_{0}))\gamma_{r}(\|e_{r}(t_{0})\|)$.
\eproof

\prooflater{Lemma \ref{lem:5-1}} For $\alpha_{s}(s)$ in $\textbf{DC}_{\textbf{2}}$,
one has $\alpha_{s}(\mu)= (\alpha_{x}(\mu))^{m}\kappa^{\frac{v_{1}}{2}}\big(\alpha_{x}(\mu)\big)\exp\big(\frac{v_{1}}{2}\sint_{0}^{\mu(t_{0})}\tau^{-2}\alpha_{x}(\tau)\mathrm{d}\tau\big)$. 
By (\ref{eq:til_e_s^i}), for $t\in\mathcal{T}_{p}$, one has
\begin{align}
\Vert\tilde{s}\Vert  \leq\varepsilon_{s}(\|\tilde{e}_{s}\|_{\mathcal{T}})(\alpha_{x}(\mu))^{-m}\kappa^{-\frac{v_{1}}{2}}\left(\alpha_{x}(\mu)\right).\label{eq:til_s-1}
\end{align}
 Note that $\varepsilon_{s}(s)=\exp\big(-\frac{v_{1}}{2}\sint_{0}^{\mu(t_{0})}\tau^{-2}\alpha_{x}(\tau)\mathrm{d}\tau\big)s$
belongs to $\mathcal{K}_{\infty}$, since $\exp\big(-\frac{v_{1}}{2}\sint_{0}^{\mu(t_{0})}\tau^{-2}\alpha_{x}(\tau)\mathrm{d}\tau\big)$
is a finite constant. Define
\begin{equation}
\tilde{r}_{1}=r_{1}-b'\otimes\varpi^{i},\label{eq:til_r_1-1}
\end{equation}
 where $b'=[1,\cdots,0]\t\in\mathbb{R}^{m-1}$ and $r_{1}$ is defined
in (\ref{eq:til_r^i_1}). Note that
\begin{equation}
e_{s}= \big [((\Phi^{-1}(\mu)\otimes I_{n})\tilde{r}_{1})\t, x_m\t\big ]\t,\label{eq:es_note}
\end{equation}
where $\Phi(\mu)$ is given in (\ref{eq:Phimu}). Taking time derivative
of $\tilde{r}_{1}$ and using (\ref{eq:r_2-1}) yield $\dot{\tilde{r}}_{1}=\alpha_{x}(\mu)(\Lambda\otimes I_{n})\tilde{r}_{1}+\delta_{x}(\mu)(A\otimes I_{n})\tilde{r}_{1}+\alpha_{x}(\mu)k_{1}(b\otimes\tilde{s})-b'\otimes\dot{\varpi}^{i}$,
where $b=[0,\cdots,1]\t\in\mathbb{R}^{m-1}$, $\Lambda$ is denoted
in (\ref{eq:Lambda}), and $A = [0, 0_{m-2}\t; 0_{m-2}, \mbox{diag}\left\{ -L_{2},\cdots,-L_{m-2}\right\} ]$. 
 Let the Lyapunov function candidate for $\tilde{r}_{1}$-dynamics
be $V_{r}(\tilde{r}_{1})=\tilde{r}_{1}\t P\tilde{r}_{1}$. Then, its
time derivative is 
$
\dot{V}_{r}(\tilde{r}_{1})=  -\alpha_{x}(\mu)\tilde{r}_{1}\t Q\tilde{r}_{1}+2\delta_{x}(\mu)\tilde{r}_{1}\t(A\t\otimes I_{n})P\tilde{r}_{1} +2k_{1}\alpha_{x}(\mu)\tilde{r}_{1}\t P(b\otimes\tilde{s})-2\tilde{r}_{1}\t P(b'\otimes\dot{\varpi}^{i})$.

Due to (\ref{eq:er_bound-1}), (\ref{eq:til_s-1}) and Young's inequality,
the terms on the right-hand side of $\dot{V}_{r}(\tilde{r}_{1})$
satisfy $-\alpha_{x}(\mu)\tilde{r}_{1}\t Q\tilde{r}_{1}\leq-v_{1}\alpha_{x}(\mu)V_{r}$,
$2\delta_{x}(\mu)\tilde{r}_{1}\t(A\t\otimes I_{n})P\tilde{r}_{1}\leq v_{3}\delta_{x}(\mu)V_{r}$,
$2k_{1}\alpha_{x}(\mu)\tilde{r}_{1}\t P(b\otimes\tilde{s})\leq V_{r}+v_{4}(\varepsilon_{s}(\|\tilde{e}_{s}\|_{\mathcal{T}}))^{2}\kappa^{-v_{1}}(\alpha_{x}(\mu))$,
and $-2\tilde{r}_{1}\t P(b'\otimes\dot{\varpi}^{i})\leq V_{r}+v_{5}\kappa^{-c^{*}}(\alpha(\mu))$,
where $v_{1}$ and is given in (\ref{eq:c^i_1}), $v_{3}=2(m-3)\lambda_{\max}(P)\lambda_{\min}^{-1}(P)$,
$v_{4}=\lambda_{\min}^{-1}(P)\lambda_{\max}^{2}(P)\vert k_{1}\vert^{2}(\alpha_{x}(\mu(t_{0})))^{-2m}$,
$v_{5}=\lambda_{\min}^{-1}(P)\lambda_{\max}^{2}(P)(\gamma_{e}(\|e_{r}(t_{0})\|))^{2}$
and we used $\Vert b\Vert=\Vert b'\Vert=1$ as well as $\|\dot{\varpi}^{i}\|\leq\|\dot{e}_{r}\|$.
Therefore, by $\textbf{DC}_{\textbf{1}}$ and the above inequalities, the
bound of $\dot{V}_{r}$ can be expressed as
\begin{align}
\dot{V}_{r} (\tilde{r}_{1})&\leq(-v_{1}\alpha_{x}(\mu)+v_{3}\delta_{x}(\mu)+2)V_{r}(\tilde{r}_{1})\nonumber \\
 & \quad+(v_{4}(\varepsilon_{s}(\|\tilde{e}_{s}\|_{\mathcal{T}}))^{2}+v_{5})\kappa^{-v_{1}}(\alpha_{x}(\mu)).\label{eq:dot_V_r-1-1}
\end{align}
Due to $\textbf{DC}_{\textbf{1}}$, 
a few facts are $\exp\big(\sint_{\tau}^{t}\delta_{x}(\mu(s))\mathrm{d}s\big)  =\alpha_{x}(\mu)(\alpha_{x}(\mu(\tau)))^{-1},\forall\tau\leq t$ and $(\alpha_{x}(\mu))^{v_{3}}\kappa^{-\frac{v_{1}}{2}}\left(\alpha_{x}(\mu)\right)\leq(\alpha_{x}(\mu(t_{0})))^{v_{3}}$. 
Then
invoking comparison lemma for (\ref{eq:dot_V_r-1-1}) yields $V_{r}(t) \leq (\lambda_{\max}(P)\|\tilde{r}_{1}(t_{0})\|^{2}+v_{4}T(\varepsilon_{s}(\|\tilde{e}_{s}\|_{\mathcal{T}}))^{2}+v_{5}T)\exp(2T)\kappa^{-\frac{v_{1}}{2}}\left(\alpha_{x}(\mu)\right)$, 
where we simply replace $V_{r}(\tilde{r}_{1}(t))$, $V_{r}(\tilde{r}_{1}(t_{0}))$
with $V_{r}(t)$, $V_{r}(t_{0})$. Therefore, utilizing the property
$\left(\ssum_{i=1}^{n}|x_{i}|\right)^{p}\leq\ssum_{i=1}^{n}|x_{i}|^{p}$
for $x_{i}\in\mathbb{R}$ where $i=1,\cdots,n$, $0<p\leq1$, the
bound of $\tilde{r}_{1}$ satisfies
\begin{equation}
\|\tilde{r}_{1}\|\leq\varepsilon_{s}^{e}(\|\tilde{e}_{s}\|_{\mathcal{T}})\kappa^{-\frac{v_{1}}{4}}\left(\alpha_{x}(\mu)\right),\label{eq:norm_til_r_1-1}
\end{equation}
 where $\varepsilon_{s}^{e}\in\mathcal{K}_{\infty}^{e}$ is $\varepsilon_{s}^{e}(s)=\exp(T)\lambda_{\min}^{-\frac{1}{2}}(P)[(v_{4}T)^{\frac{1}{2}}\varepsilon_{s}(s)+\lambda_{\max}^{\frac{1}{2}}(P)\|\tilde{r}_{1}(t_{0})\|+(v_{5}T)^{\frac{1}{2}}]$
with $\varepsilon_{s}(s)$ denoted in (\ref{eq:til_s-1}). By (\ref{eq:til_s^i-2}),  (\ref{eq:til_r^i_1}),
(\ref{eq:til_r_1-1}), and (\ref{eq:norm_til_r_1-1}), we have
 $\|(\Phi^{-1}(\mu)\otimes I_{n})\tilde{r}_{1}\|\leq\epsilon_{e}\varepsilon_{s}^{e}(\|\tilde{e}_{s}\|_{\mathcal{T}})\left(\alpha_{x}(\mu(t_{0}))\right)^{L_{m}}\kappa^{-\frac{v_{1}}{4m}}\left(\alpha_{x}(\mu)\right)$ and $\Vert x_{m}\Vert\leq\big(|k_{1}|\varepsilon_{s}(\|\tilde{e}_{s}\|_{\mathcal{T}})(\alpha_{x}(\mu(t_{0})))^{-1}+\|K\|\left(\alpha_{x}(\mu(t_{0}))\right)^{L_{m}}\big)\kappa^{-\frac{v_{1}}{4m}}\left(\alpha_{x}(\mu)\right)$,
 where $\epsilon_{e}=\max\{1,\left(\alpha_{x}(\mu(t_{0}))\right)^{-L_{2}},\cdots,\left(\alpha_{x}(\mu(t_{0}))\right)^{-L_{m}}\}$.

Summarizing  (\ref{eq:es_note}) and the above inequalities, 
the bound of $e_{s}$ satisfies (\ref{eq:norm_e_s}) with $\tilde{\varepsilon}_{s}^{e}(s)\in\mathcal{\mathcal{K}}_{\infty}^{e}$
being $\tilde{\varepsilon}_{s}^{e}(s)=\epsilon_{e}\varepsilon_{s}^{e}(s)\left(\alpha_{x}(\mu(t_{0}))\right)^{L_{m}}+|k_{1}|\varepsilon_{s}(s)(\alpha_{x}(\mu(t_{0})))^{-1}+\|K\|\left(\alpha_{x}(\mu(t_{0}))\right)^{L_{m}}$.
 \eproof

\prooflater{Proposition \ref{prop:5-1}}
According to $\textbf{DC}_{\textbf{2}}$, $\alpha_b(s)$ and $\tilde \alpha_b(s)$ can be written as 
\begin{gather}
\alpha_b(s) = |k_1^{-1}| \alpha_x(s)\exp \left( \frac{v_1}{2} \sint_{0}^s \tau ^{-2} \alpha_x(\tau) \mathrm d\tau\right),\label{eq:alpha_b-1}\\
\tilde \alpha_b(s) = |k_1| \left(\alpha_x(s)\right)^{2m-1} \exp \left( \frac{v_1}{2} \sint_{0}^s \tau ^{-2} \alpha_x(\tau) \mathrm d\tau\right).\label{eq:tilde_alpha_b}
\end{gather}
Referring to $\mathcal K_\infty$ function, a continuous function $\alpha:[0,\infty)\to [0,\infty)$ is said to belong to class $\mathcal K_\infty$ if it is strictly increasing, $\alpha(0)=0$, and $\alpha(r)\to \infty$ as $r\to \infty$. 
Given that $\alpha_x(s)$ is first-order differentiable, it follows that both $ \alpha_b(s)$ and $\tilde \alpha_b(s)$ in \eqref{eq:alpha_b-1} and \eqref{eq:tilde_alpha_b} are continuous for $s\in[0,\infty)$. For $\alpha_b(s)$ in \eqref{eq:alpha_b-1}, we have $\alpha_b(0)=0$ and  $\alpha_b(s)\geq |k_1^{-1}|\alpha_x(s)$ since $\exp \left( \frac{v_1}{2} \sint_{0}^s \tau ^{-2} \alpha_x(\tau) \mathrm d\tau\right)\geq 1$. Therefore, as $\alpha_x(s) \in \mathcal{K}_\infty$, we conclude $\lim_{s \to \infty} \alpha_b(s) = \infty$.
Differentiating $\alpha_b(s)$ with respect to $s$ yields $\frac{\partial \alpha_b(s)}{\partial s } =  |k_1^{-1}| \frac{\partial \alpha_x(s)}{\partial s}\exp \left( \frac{v_1}{2} \sint_{0}^s \tau ^{-2} \alpha_x(\tau) \mathrm d\tau\right) + |k_1^{-1}| \left(\alpha_x(s)\right)^2\exp \left( \frac{v_1}{2} \sint_{0}^s \tau ^{-2} \alpha_x(\tau) \mathrm d\tau\right)\frac{v_1}{2} s^{-2}$. 

Using the fact that $\frac{\partial \alpha_x(s)}{\partial s}>0$ leads to $\frac{\partial \alpha_b(s)}{\partial s }>0$ for $s\in\mathbb R_{\geq 0}$, which implies $\alpha_b(s)$ is strictly increasing with respect to $s$. Combining this with $\alpha_b(0) = 0$ and $\lim_{s \to \infty} \alpha_b(s) = \infty$, we establish that $\alpha_b(s) \in \mathcal{K}_\infty$. The same methodology can be applied to show that $\tilde \alpha_b(s)$ belongs to class $\mathcal{K}_\infty$. This concludes the proof.
\eproof 

\prooflater{Lemma \ref{lem:5-2}} Let the Lyapunov function candidate
for $\tilde{e}_{s}$-dynamics in (\ref{eq:dot_omega-2}) be
$
V_{s}(\tilde{e}_{s})=\tilde{e}_{s}\t\tilde{e}_{s}/2$. 
According to Theorem \ref{the:auxi_genera}, $\dot{\varpi}^{i}$ satisfies
$\|\dot{\varpi}^{i}\|\leq\|\dot{\varpi}\|\leq\|\dot{e}_{r}\|\leq\gamma_{\zeta}(\mu)\|e_{r}\|$
with $\gamma_{\zeta}(s)$ denoted in (\ref{eq:alpha_r}). Utilizing
(\ref{eq:psi^i}), $\textbf{DC}_{\textbf{2}}$, (\ref{eq:alpha_b}) and Young's
inequality, we have $\tilde{e}_{s}\t\alpha_{s}(\mu)B(\mu)\varphi(x,d)\leq\alpha_{b}(\mu)((\psi^{2}(x)\Vert\tilde{e}_{s}\Vert^{2}+h^{2}(\|d\|)/4)$
and $\tilde{e}_{s}\t\alpha_{s}(\mu)\dot{\varpi}^{i}\leq\alpha_{b}(\mu)\Vert\tilde{e}_{s}\Vert^{2}+\tilde{\alpha}_{b}(\mu)(\gamma_{\zeta}(\mu))^{2}\Vert e_{r}\Vert^{2}/4$.
Then the time-derivative of $V_{s}(\tilde{e}_{s})$ along $\tilde{e}_{s}$-dynamics
(\ref{eq:dot_omega-2}) with the controller (\ref{eq:u^i}) satisfies 
$
\dot{V}_{s}(\tilde{e}_{s})
  \leq-v\alpha_{b}(\mu)\|\tilde{e}_{s}\|^{2}+\alpha_{b}(\mu)(h(\|d\|))^{2}/4 +\tilde{\alpha}_{b}(\mu)(\gamma_{\zeta}(\mu))^{2}\Vert e_{r}\Vert^{2}/4$.
  
Therefore, 
the closed-loop system of $\tilde{e}_{s}$ admit a prescribed-time
ISS Lyapunov function in the form of $\textbf{C}_{\textbf{2}}$ with
the prescribed-time convergent gain, prescribed-time ISS gain and
ISS gain given as in (\ref{eq:con_1}). 

What remains is to prove that control input $u$ in (\ref{eq:u^i})
satisfies $\textbf{C}_{\textbf{3}}'$. By
(\ref{eq:norm_til_r_1-1}), $\|\pi\|\leq\varepsilon_{1}+\varepsilon_{2}\|\tilde{e}_{s}\|_{\mathcal{T}}+\varepsilon_{3}(\alpha_{x}(\mu))^{m}\|e_{r}\| $, 
 where $\varepsilon_{1}$, $\varepsilon_{2}$ and $\varepsilon_{3}$
are  finite constants. And $\|B^{-1}(\mu)\delta_{s}(\mu)\tilde{s}\|
\leq\varepsilon_{4}\|\tilde{e}_{s}\|_{\mathcal{T}}$
for some finite constant $\varepsilon_{4}$. Since $\psi(x)$ is bounded
for $x\in\mathbb{R}^{mn}$, one has $-(v+\psi^{2}(x)+1)\mbox{sign}(k_{1})\tilde{e}_{s}\leq\varepsilon_{5}\|\tilde{e}_{s}\|_{\mathcal{T}}$
for some positive finite constant $\varepsilon_{5}$. Therefore, the controller $u$ in (\ref{eq:u^i})
satisfies $\textbf{C}_{\textbf{3}}'$ with $\tilde{\gamma}_{s}(s)$
and $\tilde{\gamma}_{u}$ in (\ref{eq:til_gamma_s-1}) by letting
$\varepsilon_{1}'=\varepsilon_{1}$, $\varepsilon_{2}'=\varepsilon_{2}+\varepsilon_{4}+\varepsilon_{5}$
and $\varepsilon_{3}'=\varepsilon_{3}$.\eproof

\prooflater{Theorem \ref{thm:5-1}} In order to
invoke Theorem \ref{thm:2}, we must check the
conditions (\ref{eq:con2-1}) and (\ref{eq:con3-3}).
Due to Lemma \ref{lem:5-2}, one obtains  $\tilde{\sigma}_{\tilde{d}}(s)=\alpha_{b}(s)/4$.
Therefore,
\begin{align*}
\sup_{s\in\mathbb{R}_{\geq0}}[\tilde{\sigma}_{\tilde{d}}(s)/(\exp(c)\tilde{\alpha}_{\tilde{s}}(s))]=(8v\exp(c))^{-1}<\infty,
\end{align*} 
which means (\ref{eq:con2-1}) is satisfied. 
By Theorem \ref{the:auxi_genera}
and $\textbf{DC}_{\textbf{1}}$, one has $\varrho_{r}(c,s)=(2c/c_{1})^{\frac{1}{2}}\exp\left(-v_{1}\sint_{0}^{s}\tau^{-2}\alpha_{x}(\tau)\mathrm{d}\tau\right)$, 
 where $c_{1}$, $c_{2}$ and $c^{*}$ are denoted in (\ref{eq:c_1}). Therefore, $\sup_{s\in\mathbb{R}_{\geq0}}\left[\tilde{\sigma}_{\tilde{r}}(s)/(\tilde{\alpha}_{\tilde{s}}(s)(\sigma_{\tilde{r}}(\varrho_{r}(c,s)))^{-1})\right]=\sup_{s\in\mathbb{R}_{\geq0}}\left[\varepsilon(\alpha_{x}(s))^{2m}\exp\left(-2v_{1}\sint_{0}^{s}\tau^{-2}\alpha_{x}(\tau)\mathrm{d}\tau\right)\right] \leq\varepsilon$,
where $\varepsilon=ck_{1}^{2}v_{1}^{2}(2\lambda_{N}+\varrho_{c})^{2}/(4vc_{1}(c^{*})^{2})$
and we used $m/v_{2}\leq1$. Consequently, condition (\ref{eq:con3-3}) is satisfied.
By Theorem \ref{thm:2}, $\tilde{e}_{s}$ is bounded. In Theorem \ref{the:auxi_genera},
we have proved $\gamma_{\zeta}(s)=\mathcal{S}[\varrho_{r}(c,s)]$. 

Furthermore, $  \sup_{s\in\mathbb{R}_{\geq0}}[\tilde{\gamma}_{u}/(\varrho_{r}(c,s))^{-1}] =\sup_{s\in\mathbb{R}_{\geq0}}\varepsilon_{3}'(2c/c_{1})^{\frac{1}{2}}(\alpha_{x}(s))^{m}\exp\left(-v_{1}\sint_{0}^{s}\tau^{-2}\alpha_{x}(\tau)\mathrm{d}\tau\right) \leq\varepsilon_{3}'(2c/c_{1})^{\frac{1}{2}}$.
By Lemmas \ref{lem:5-1} and \ref{lem:5-2}, all conditions in Theorem
\ref{thm:2} are satisfied. Invoking Theorem \ref{thm:2}
completes the proof. \eproof

\prooflater{Lemma \ref{lem:6-1}} Due to (\ref{eq:til_x}), one has 
\[
e_{s}=\left[\begin{array}{c}
x_{1}-\varpi^{i}\\
x_{2}\\
\vdots\\
x_{m}\\
\hat{\theta}\\
\xi_{f}
\end{array}\right]=\left[\begin{array}{c}
\tilde{x}_{1}\\
\tilde{x}_{2}+\tilde{\xi}_{2}+\xi_{1}\\
\vdots\\
\tilde{x}_{m}+\tilde{\xi}_{m}+\xi_{m-1}\\
\hat{\theta}\\
\xi_{f}
\end{array}\right].
\]
By (\ref{eq:omega eta}), then
\begin{equation}\label{eq:til_x_1}
	\Vert\tilde{x}_{q}\Vert  \leq\|\tilde{e}_{s}\|_{\mathcal{T}}(\alpha_{\xi}(\mu))^{-L_{q}},\;  \Vert\tilde{\xi}_{q}\Vert  \leq\|\tilde{e}_{s}\|_{\mathcal{T}}(\alpha_{\xi}(\mu))^{-L_{q}}
\end{equation}
where we used $\|\tilde{x}_{q}\|\leq\|\omega\|(\alpha_{\xi}(\mu))^{-L_{q}}$,
$\|\tilde{\xi}_{q}\|\leq\|\eta\|(\alpha_{\xi}(\mu))^{-L_{q}}$, $\|\omega\|\leq\|\tilde{e}_{s}\|\leq\|\tilde{e}_{s}\|_{\mathcal{T}}$
and $\|\eta\|\leq\|\tilde{e}_{s}\|\leq\|\tilde{e}_{s}\|_{\mathcal{T}}$.
Due to (\ref{eq:til_x}), the bounds of $\xi_{1}$ and $\xi_{2f}$
satisfy $	\|\xi_{1}\| \leq c_{1}\|\tilde{e}_{s}\|_{\mathcal{T}}(\alpha_{\mu}(\mu))^{-L_{2}}$ and $\|\xi_{2f}\|\leq\Vert\tilde{\xi}_{2}\Vert+\Vert\xi_{1}\Vert\leq(1+c_{1})\|\tilde{e}_{s}\|_{\mathcal{T}}(\alpha_{\mu}(\mu))^{-L_{2}}$.
As a result, the bound of $x_{2}$ satisfies
\begin{align}
\Vert x_{2}\Vert &  \leq(2+c_{1})\|\tilde{e}_{s}\|_{\mathcal{T}}(\alpha_{\xi}(\mu))^{-L_{2}}.\label{eq:x_2}
\end{align}
By Assumption \ref{ass:varphi-1}, one has $\|\hat{\theta}\varphi_{2}(x_{2})\|\leq|\hat{\theta}|\bar{\psi}_{2}\|x_{2}\|\leq\bar{\psi}_{2}(\|\tilde{e}_{s}\|_{\mathcal{T}}+|\theta|)\|x_{2}||$,
where we used $|\hat{\theta}|\leq|\tilde{\theta}|+|\theta|$ and $|\tilde{\theta}|\leq\|\tilde{e}_{s}\|_{\mathcal{T}}$.
Then, the bound of $\xi_{2}$ and $\xi_{3f}$ satisfy $\Vert\xi_{2}\Vert  \leq(c_{2}+\bar{\psi}_{2}|\theta|(2+c_{1})+v_{2})\|\tilde{e}_{s}\|_{\mathcal{T}}(\alpha_{\xi}(\mu))^{-L_{3}}+\bar{\psi}_{2}(2+c_{1})\|\tilde{e}_{s}\|_{\mathcal{T}}^{2}(\alpha_{\xi}(\mu))^{-L_{3}}$ and $\Vert\xi_{3f}\Vert  \leq(c_{2}+\bar{\psi}_{2}|\theta|(2+c_{1})+v_{2}+1)\|\tilde{e}_{s}\|_{\mathcal{T}}(\alpha_{\xi}(\mu))^{-L_{3}}+\bar{\psi}_{2}(2+c_{1})\|\tilde{e}_{s}\|_{\mathcal{T}}^{2}(\alpha_{\xi}(\mu))^{-L_{3}}$.
Similarly, we can derive that
\begin{equation}
\begin{gathered}\Vert x_{q}\Vert\leq\varepsilon_{xq}(\|\tilde{e}_{s}\|_{\mathcal{T}})(\alpha_{\xi}(\mu))^{-L_{q}},\;q=3,\cdots,m,\\
\Vert\xi_{q}\Vert\leq\varepsilon_{\xi q}(\|\tilde{e}_{s}\|_{\mathcal{T}})(\alpha_{\xi}(\mu))^{-L_{q}},\;q=3,\cdots,m,\\
\Vert\xi_{qf}\Vert\leq\varepsilon_{\xi q}'(\|\tilde{e}_{s}\|_{\mathcal{T}})(\alpha_{\xi}(\mu))^{-L_{q}},\;q=4,\cdots,m,
\end{gathered}
\label{eq:conver_x_xi-1}
\end{equation}
where $\varepsilon_{xq}$, $\varepsilon_{\xi q}$ and $\varepsilon_{\xi q}'$
are some $\mathcal{K}_{\infty}$ functions.
By (\ref{eq:Hurwitz}), (\ref{eq:omega eta}), (\ref{eq:x_2}) and
(\ref{eq:conver_x_xi-1}), $\tau$ in (\ref{eq:tau_q^i}) satisfies
\begin{equation}
\vert\tau\vert\leq\ssum_{j=2}^{m}\Vert\omega_{j}\Vert(\alpha_{\xi}(\mu))^{L_{j}}\bar{\psi}_{j}\Vert x_{j}\Vert\leq\varepsilon_{\tau}(\|\tilde{e}_{s}\|_{\mathcal{T}}),\label{eq:norm_tau-1}
\end{equation}
where $\varepsilon_{\tau}\in\mathcal{K_{\infty}}$. Define the Lyapunov
function candidate for $\hat{\theta}$-dynamics as $U=\frac{1}{2}\hat{\theta}^{2}$.
Its time-derivative along $\hat{\theta}$-dynamics is 
$
\dot{U}\leq-(2\sigma-1)\alpha_{\xi}(\mu)U+(\alpha_{\xi}(\mu))^{-1}(\varepsilon_{\tau}(\|\tilde{e}_{s}\|_{\mathcal{T}}))^{2}/2
$. 
Let $U_{s}=(\alpha_{\xi}(\mu))^{2}U$ and
\begin{equation}
\sigma=(3+\sigma')/2\label{eq:sigma}
\end{equation}
with any $\sigma'>0$. By $U$-dynamics and $\textbf{DC}_{\xi}$,
$\dot{U}_{s}\leq-\sigma'\alpha_{\xi}(\mu)U_{s}+\frac{1}{2}\alpha_{\xi}(\mu)(\varepsilon_{\tau}(\|\tilde{e}_{s}\|_{\mathcal{T}}))^{2}$.
Invoking comparison lemma, one has $U_{s}(t)\leq\kappa^{-\sigma'}(\alpha_{\xi}(\mu))U_{s}(t_{0})+\frac{1}{2\sigma'}(\varepsilon_{\tau}(\|\tilde{e}_{s}\|_{\mathcal{T}}))^{2}$,
then $\hat{\theta}$ satisfies $\vert\hat{\theta}(t)\vert  \leq\big(\alpha_{\xi}(\mu(t_{0}))|\hat{\theta}(t_{0})|+\left(2\sigma'\right)^{-\frac{1}{2}}\varepsilon_{\tau}(\|\tilde{e}_{s}\|_{\mathcal{T}})\big)(\alpha_{\xi}(\mu))^{-1} \leq\varepsilon_{\hat{\theta}}^{e}(\|\tilde{e}_{s}\|_{\mathcal{T}})(\alpha_{\xi}(\mu))^{-1}$
for $\varepsilon_{\hat{\theta}}^{e}\in\mathcal{K}_{\infty}^{e}$.
Therefore, (\ref{eq:til_x_1})-(\ref{eq:conver_x_xi-1}) lead to
(\ref{eq:norm_e_s-2}). \eproof

\prooflater{Theorem \ref{thm:6-2}} Let $V_{i}=\sum_{q=1}^{i}\frac{1}{2}\left(\omega_{q}\t\omega_{q}+\eta_{q+1}\t\eta_{q+1}\right)$.
We use the backstepping technique to prove
\begin{align}
\dot{V}_{i}\leq & -\ssum_{j=1}^{i}\bar{c}_{j}\alpha_{\xi}(\mu)\Vert\omega_{j}\Vert^{2}-\ssum_{j=1}^{i}\bar{\upsilon}_{j+1}\alpha_{\xi}(\mu)\Vert\eta_{i+1}\Vert^{2}\nonumber \\
 & +\alpha_{\xi}(\mu)\Vert\omega_{i+1}\Vert^{2}/2+\tilde{\theta}\ssum_{j=2}^{i}\tau_{j}+\alpha_{\xi}(\mu)\ssum_{j=1}^{i}\pi_{j}\nonumber \\
 & +(\alpha_{\xi}(\mu))^{2L_{2}+1}\Vert\dot{\varpi}^{i}\Vert^{2}/2\label{eq:Vq}
\end{align}
holds for $i=1,\cdots,m-1$, where $\bar{c}_{j}$ and 
$\bar{\upsilon}_{j+1}$ are some positive finite constants, $\tau_{j}$
is denoted in (\ref{eq:tau_q^i}) and $\pi_{j}$ will be defined later. 

{\emph{Step 1}}: The derivative of $\omega_{1}\t\omega_{1}$
is $\omega_{1}\t\dot{\omega}_{1}=\omega_{1}\t(\alpha_{\xi}(\mu))^{L_{1}}\big(\xi_{1}+\tilde{x}_{2}+\tilde{\xi}_{2}-\dot{\varpi}^{i}+L_{1}\delta_{\xi}(\mu)\tilde{x}_{1}\big)$, 
where $\delta_{\xi}(\mu)=\frac{\mathrm{d}\alpha_{\xi}(\mu)}{\mathrm{d}\mu}\mu^{2}(\alpha_{\xi}(\mu))^{-1}$. 
By Young's inequality and $\textbf{DC}_{\xi}$, we have $	\omega_{1}\t(\alpha_{\xi}(\mu))^{L_{1}}L_{1}\alpha_{\xi}'(\mu)\tilde{x}_{1}\leq L_{1}\alpha_{\xi}(\mu)\Vert\omega_{1}\Vert^{2}$, $	\omega_{1}\t(\alpha_{\xi}(\mu))^{L_{1}}\tilde{x}_{2}\leq\frac{1}{2}\alpha_{\xi}(\mu)(\Vert\omega_{1}\Vert^{2}+\Vert\omega_{2}\Vert^{2})$, $	 \omega_{1}\t(\alpha_{\xi}(\mu))^{L_{1}}\tilde{\xi}_{2}\leq\frac{1}{2}\alpha_{\xi}(\mu)(\Vert\omega_{1}\Vert^{2}+\Vert\eta_{2}\Vert^{2})$ and $\omega_{1}\t(\alpha_{\xi}(\mu))^{L_{1}}\dot{\varpi}^{i}\leq\frac{1}{2}\alpha_{\xi}(\mu)\Vert\omega_{1}\Vert^{2}+\frac{1}{2}(\alpha_{\xi}(\mu))^{2L_{2}+1}\Vert\dot{\varpi}^{i}\Vert^{2}$. 
Let $c_{1}=\bar{c}_{1}+L_{1}+\frac{3}{2}$ with any $\bar{c}_{1}\geq\frac{1}{2}\sigma$,
where $\sigma$ is denoted in (\ref{eq:sigma}). Substituting $\xi_{1}$
in (\ref{eq:til_x}) and the above inequalities into $\omega_{1}\t\dot{\omega}_{1}$
yields $\omega_{1}\t\dot{\omega}_{1}\leq  -\bar{c}_{1}\alpha_{\xi}(\mu)\Vert\omega_{1}\Vert^{2}+\alpha_{\xi}(\mu)(\Vert\omega_{2}\Vert^{2}+\Vert\eta_{2}\Vert^{2})/2 +(\alpha_{\xi}(\mu))^{2L_{2}+1}\Vert\dot{\varpi}^{i}\Vert^{2}/2$. 

Let $\upsilon_{2}=\bar{\upsilon}_{2}+L_{2}+\rho_{2}+\frac{1}{2}$
with $\bar{\upsilon}_{2}\geq\frac{1}{2}\sigma$ and $\rho_{2}$ sufficiently
large. The term $\eta_{2}\t\dot{\eta}_{2}$ satisfies $\eta_{2}\t\dot{\eta}_{2} 
 \leq-(\bar{\upsilon}_{2}+1/2)\alpha_{\xi}(\mu)\Vert\eta_{2}\Vert^{2}+\alpha_{\xi}(\mu)\pi_{1}$, 
where
\begin{equation}
\pi_{1}=\Vert\left(\alpha_{\xi}(\mu)\right)^{L_{3}}\dot{\xi}_{1}\Vert^{2}/(2\rho_{2}).\label{eq:pi1}
\end{equation}
Note that $\dot{\xi}_{1}$ exists for $t\in\mathcal{T}_{p}$ and can
be expressed as
$
\dot{\xi}_{1}=-c_{1}\delta_{\xi}(\mu)\alpha_{\xi}(\mu)\tilde{x}_{1}-c_{1}\alpha_{\xi}(\mu)(x_{2}-\dot{\varpi}^{i})
$.
Consider the Lyapunov function candidate as $V_{1}=\frac{1}{2}\left(\Vert\omega\Vert^{2}+\Vert\eta_{2}\Vert^{2}\right)$.
 Then $\dot{V}_{1}$
satisfies
\begin{align}
\dot{V}_{1}\leq & -\bar{c}_{1}\alpha_{\xi}(\mu)\Vert\omega_{1}\Vert^{2}-\bar{\upsilon}_{2}\alpha_{\xi}(\mu)\Vert\eta_{2}\Vert^{2}+\alpha_{\xi}(\mu)\Vert\omega_{2}\Vert^{2}/2\nonumber \\
 & +(\alpha_{\xi}(\mu))^{2L_{2}+1}\Vert\dot{\varpi}^{i}\Vert^{2}+\alpha_{\xi}(\mu)\pi_{1}/2.\label{eq:dot_V_1-1}
\end{align}

{\emph{Step q}}: Suppose (\ref{eq:Vq}) holds for $i=q-1$.
Let $c_{q}=\bar{c}_{q}+L_{q}+2$ with any $\bar{c}_{q}\geq\frac{1}{2}\sigma$.
Substituting $\xi_{q}$ in (\ref{eq:til_x}) into $\omega_{q}\t\dot{\omega}_{q}$
yields
$ 
\omega_{q}\t\dot{\omega}_{q} 
 \leq-\left(\bar{c}_{q}+1/2\right)\alpha_{\xi}(\mu)\Vert\omega_{q}\Vert^{2}+\tilde{\theta}\tau_{q} \quad+\alpha_{\xi}(\mu)(\Vert\omega_{q+1}\Vert^{2}+\Vert\eta_{q+1}\Vert^{2})/2$.
Let $\upsilon_{q+1}=\bar{\upsilon}_{q+1}+L_{q+1}+\rho_{q+1}+\frac{1}{2}$
with $\bar{\upsilon}_{q+1}\geq\frac{1}{2}\sigma$ and $\rho_{q+1}$
sufficiently large. Then, $\eta_{q+1}\t\dot{\eta}_{q+1}$ satisfies
$
\eta_{q+1}\t\dot{\eta}_{q+1}\leq-(\bar{\upsilon}_{q+1}+1/2)\alpha_{\xi}(\mu)\Vert\eta_{q+1}\Vert^{2}+\alpha_{\xi}(\mu)\pi_{q}$,
where $\pi_{q}=\Vert\left(\alpha_{\xi}(\mu)\right)^{L_{q+2}}\dot{\xi}_{q}\Vert^{2}/(2\rho_{q+1})$.

By Assumption \ref{ass:varphi-1}, $\dot{\xi}_{q}$ can be expressed
as  $\dot{\xi}_{q}=  -c_{q}\frac{\mathrm{d}\alpha_{\xi}(\mu)}{\mathrm{d}\mu}\mu^{2}\tilde{x}_{q}-c_{q}\alpha_{\xi}(\mu)\dot{\tilde{x}}_{q}-\dot{\hat{\theta}}\psi_{q}(x_{q})x_{q} -\hat{\theta}\psi_{q}(x_{q})\dot{x}_{q}-\hat{\theta}\frac{\mathrm{d}\psi_{q}(x_{q}(t))}{\mathrm{d}t}x_{q}-\upsilon_{q}\frac{\mathrm{d}\alpha_{\xi}(\mu)}{\mathrm{d}\mu}\mu^{2}\tilde{\xi}_{q} -\upsilon_{q}\alpha_{\xi}(\mu)\dot{\tilde{\xi}}_{q}$. 
Note that $V_{q}=V_{q-1}+\frac{1}{2}\left(\omega_{q}\t\omega_{q}+\eta_{q+1}\t\eta_{q+1}\right)$.
By (\ref{eq:dot_V_1-1}),
$\dot{V}_{q}$ can be expressed as $\dot{V}_{q}\leq  -\ssum_{j=1}^{q}\bar{c}_{j}\alpha_{\xi}(\mu)\Vert\omega_{j}\Vert^{2}-\ssum_{j=1}^{q}\bar{\upsilon}_{j+1}\alpha_{\xi}(\mu)\Vert\eta_{q+1}\Vert^{2} +\alpha_{\xi}(\mu)\Vert\omega_{q+1}\Vert^{2}/2+\tilde{\theta}\ssum_{j=2}^{q}\tau_{j}+\alpha_{\xi}(\mu)\ssum_{j=1}^{q}\pi_{j} +(\alpha_{\xi}(\mu))^{2L_{2}+1}\Vert\dot{\varpi}^{i}\Vert^{2}/2$.
Thus, the claim (\ref{eq:Vq}) holds for $i=1,\cdots,m-1$.

{\emph {Step m}}: Let Lyapunov function candidate be $V_{m}=V_{m-1}+\frac{1}{2}(\omega_{m}\t\omega_{m}+\tilde{\theta}^{2})$.
Substituting (\ref{eq:u^i-1}) and $\dot{\hat{\theta}}$ in (\ref{eq:adaptive_law})
into $\dot{V}_{m}$ yields $\dot{V}_{m}\leq  -\alpha_{\xi}(\mu)\ssum_{j=1}^{m}\bar{c}_{j}\Vert\omega_{j}\Vert^{2}-\alpha_{\xi}(\mu)\ssum_{j=1}^{m-1}\bar{\upsilon}_{j+1}\Vert\eta_{q+1}\Vert^{2} +\sigma\alpha_{\xi}(\mu)\tilde{\theta}\hat{\theta}+\alpha_{\xi}(\mu)\ssum_{j=1}^{m-1}\pi_{j} +(\alpha_{\xi}(\mu))^{2L_{2}+1}\Vert\dot{\varpi}^{i}\Vert^{2}/2$,
where we used (\ref{eq:tau_q^i}).
For a given $h>0$, define $\Omega(h)=\left\{ \tilde{e}_{s}\in\mathbb{R}^{2mn+1-n}\mid\|\tilde{e}_{s}\|^{2}\leq h^{2}\right\} $.
Within $\Omega(h)$, one has $\|\tilde{e}_{s}\|\leq h$, $\|\omega_{j}\|_{\mathcal{}}\leq h$,
$\|\eta_{j}\|\leq h$ and $|\tilde{\theta}|\leq h$. Let constants
$c_{j}$ for $j=1,\cdots,m$ and $\upsilon_{j}$ for $j=2,\cdots,m$
be defined in the proof of Lemma \ref{lem:6-1}. Define constant vectors
$\tilde{c}_{j}=[c_{1},\cdots,c_{j}]\t$ for $j=1,\cdots,m$, $\tilde{\upsilon}_{j}=[\upsilon_{2},\cdots,\upsilon_{j}]\t$
for $j=2,\cdots,m$ and
\begin{equation}
\tilde{c}_{s,j}=[\tilde{c}_{j}\t,\tilde{\upsilon}_{j+1}\t,\sigma,\theta,h]\t\label{eq:c_s^j}
\end{equation}
for $j=2,\cdots,m$. Following the proof of Lemma
\ref{lem:6-1}, one has 
$  \|\tilde{x}_{1}\|\leq\epsilon_{1}(h)(\alpha_{\xi}(\mu))^{-L_{1}}$, $ \|x_{2}\|\leq\epsilon_{2}(\tilde{c}_{1},h)(\alpha_{\xi}(\mu))^{-L_{2}}$, $  \Vert x_{q}\Vert\leq\epsilon_{q}(\tilde{c}_{q-1},\tilde{\upsilon}_{q-1},h)(\alpha_{\xi}(\mu))^{-L_{q}}$, $ \|\xi_{1}\|\leq\bar{\xi}_{1}(\tilde{c}_{1},h)(\alpha_{\xi}(\mu))^{-L_{1}}$, $  \Vert\xi_{q}\Vert\leq\bar{\xi}_{q}(\tilde{c}_{q-1},\tilde{\upsilon}_{q-1},h,\theta)(\alpha_{\xi}(\mu))^{-L_{q}}$, 
$\|\xi_{2f}\|\leq\bar{\xi}_{2f}(\tilde{c}_{1},h)(\alpha_{\xi}(\mu))^{-L_{2}}$, $\Vert\xi_{qf}\Vert\leq\bar{\xi}_{qf}(\tilde{c}_{q-1},\tilde{\upsilon}_{q-1},h,\theta)(\alpha_{\xi}(\mu))^{-L_{q}}$ and $|\hat{\theta}|\leq\Theta(\theta,\sigma,h)(\alpha_{\xi}(\mu))^{-1}$
with some functions $\epsilon_{q}$ for $q=1,\cdots,m$, $\bar{\xi}_{q}$,
$\bar{\xi}_{qf}$ for $q=2,\cdots,m$ and $\Theta$. From Theorem
\ref{the:auxi_genera}, one has
\begin{equation}
\|\dot{\varpi}^{i}\|\leq\|\dot{e}_{r}\|\leq\gamma_{\zeta}(\mu)\|e_{r}\|\label{eq:dot_varpi-1}
\end{equation}
with $\gamma_{\zeta}(\mu)=(2\lambda_{N}+\varrho_{c})\alpha(\mu)$.
According to  (\ref{eq:dot_varpi-1}),
$\dot{\xi}_{1}$ satisfies $\|\dot{\xi}_{1}\|\leq c_{1}(\epsilon_{1}(h)+\epsilon_{2}(\tilde{c}_{1},h))(\alpha_{\xi}(\mu))^{-L_{3}}+c_{1}\alpha_{\xi}(\mu)\gamma_{\zeta}(\mu)\|e_{r}\|$.
Thus, $\pi_{1}$ in (\ref{eq:pi1}) satisfies 
\begin{equation}
\pi_{1}\leq\frac{\Xi_{1}(\tilde{c}_{s,1})}{2\rho_{2}}+\frac{\epsilon_{r,1}(\tilde{c}_{s,1})}{2\rho_{2}}(\alpha_{\xi}(\mu))^{2L_{2}}(\alpha(\mu))^{2}\|e_{r}\|^{2},\label{eq:pi_1}
\end{equation}
where $\Xi_{1}(\tilde{c}_{s,1})=c_{1}^{2}h^{2}(3+c_{1})^{2}$ and
$\varepsilon_{r,1}(\tilde{c}_{s,1})=c_{1}^{2}(2\lambda_{N}+\varrho_{c})^{2}$.
According to the above inequalities  and Remark \ref{rem:6-1},
the terms on the right hand side of $\dot{\xi}_{2}$ 
satisfy $-c_{2}\frac{\mathrm{d}\alpha_{\xi}(\mu)}{\mathrm{d}\mu}\mu^{2}\tilde{x}_{2}\leq c_{2}h(\alpha_{\xi}(\mu))^{-L_{4}}$, $-c_{2}\alpha_{\xi}(\mu)\dot{\tilde{x}}_{2}\leq c_{2}(\epsilon_{3}+\vert\theta\vert\bar{\psi}_{2}\epsilon_{2}\alpha_{\xi}(\mu(t_{0}))+\upsilon_{2}h)(\alpha_{\xi}(\mu))^{-L_{4}}$, $-\dot{\hat{\theta}}\psi_{2}(x_{2})x_{2}\leq(\varepsilon_{\tau}(h)+\sigma\Theta)\bar{\psi}_{2}\epsilon_{2}(\alpha_{\xi}(\mu))^{-L_{2}}$, $-\hat{\theta}\psi_{2}(x_{2})\dot{x}_{2}\leq\Theta\bar{\psi}_{2}(\epsilon_{3}+\alpha_{\xi}(\mu(t_{0}))\vert\theta\vert\bar{\psi}_{2}\epsilon_{2}) (\alpha_{\xi}(\mu))^{-L_{2}}$, $-\hat{\theta}\frac{\mathrm{d}\psi_{2}(x_{2}(t))}{\mathrm{d}t}x_{2}\leq\text{\ensuremath{\Theta}}\tilde{\psi}_{2}\epsilon_{2}(\alpha_{\xi}(\mu))^{-L_{1}}$, $-\upsilon_{2}\frac{\mathrm{d}\alpha_{\xi}(\mu)}{\mathrm{d}\mu}\mu^{2}\tilde{\xi}_{2}\leq\upsilon_{2}h(\alpha_{\xi}(\mu))^{-L_{4}}$ and $-\upsilon_{2}\alpha_{\xi}(\mu)\dot{\tilde{\xi}}_{2}\leq(\alpha_{\xi}(\mu)^{-L_{4}}\upsilon_{2}^{2}h+\upsilon_{2}\alpha_{\xi}(\mu)\|\dot{\xi}_{1}\|$,
where $\varepsilon_{\tau}$ is given in (\ref{eq:norm_tau-1}). As
a result, upon utilizing the above inequalities and (\ref{eq:pi_1}),
$\pi_{2}$  satisfies
\begin{equation}
\pi_{2}\leq\frac{\Xi_{2}(\tilde{c}_{s,2})}{\rho_{3}}+\frac{\varepsilon_{r,2}(\tilde{c}_{s,2})}{\rho_{3}}(\alpha_{\xi}(\mu))^{2L_{2}}(\alpha(\mu))^{2}\|e_{r}\|^{2}, \label{eq:pi_2}
\end{equation}
where $\tilde{c}_{s,2}$ is denoted in (\ref{eq:c_s^j}), $\Xi_{2}(\tilde{c}_{s,2})$
and $\varepsilon_{r,2}(\tilde{c}_{s,2})$ are independent of $\upsilon_{3}$.
Similarly, by utilizing the above inequalities,
(\ref{eq:pi_1}) and (\ref{eq:pi_2}), step by step, it can be derived
that
\begin{equation}
\pi_{j}\leq\frac{\Xi_{j}(\tilde{c}_{s,j})}{\rho_{j+1}}+\frac{\varepsilon_{r,j}(\tilde{c}_{s,j})}{\rho_{j+1}}(\alpha_{\xi}(\mu))^{2L_{2}}(\alpha(\mu))^{2}\|e_{r}\|^{2}\label{eq:pi_j}
\end{equation}
for $j=3,\cdots,m-1$, where $\Xi_{j}(\tilde{c}_{s,j})$ and $\varepsilon_{r,j}(\tilde{c}_{s,j})$
are independent of $\upsilon_{j+1}$. By (\ref{eq:dot_varpi-1}), 
$
(\alpha_{\xi}(\mu))^{2L_{2}+1}\Vert\dot{\varpi}^{i}\Vert^{2}/2\leq\tilde{\varepsilon}_{r}(\alpha_{\xi}(\mu))^{2L_{2}+1}(\alpha(\mu))^{2}\Vert e_{r}\Vert^{2}$
for some $0<\tilde{\varepsilon}_{r}<\infty$. Note that
$
V_{m}=\frac{1}{2}\big(\ssum_{j=1}^{m}\omega_{j}\t\omega_{j}+\ssum_{j=2}^{m}\eta_{j}\t\eta_{j}+\tilde{\theta}^{2}\big)=\frac{1}{2}\tilde{e}_{s}\t\tilde{e_{s}}
$.
Using (\ref{eq:pi_1}), (\ref{eq:pi_2}), and  (\ref{eq:pi_j}),
the bound of $\dot{V}_{m}$  satisfies
$
\dot{V}_{m}\leq  -\iota_{1}\alpha_{\xi}(\mu)V_{m}+\alpha_{\xi}(\mu)(\sigma\theta^{2}/2+\ssum_{j=1}^{m-1}\Xi_{j}/\rho_{j+1}) +\iota_{2}(\alpha_{\xi}(\mu))^{2L_{2}+1}(\alpha(\mu))^{2}\Vert e_{r}\text{\ensuremath{\Vert^{2}}}$, 
where $\iota_{1}=\min\{2\bar{c}_{1},\cdots,2\bar{c}_{m},2\bar{\upsilon}_{2},\cdots,2\bar{\upsilon}_{m},\sigma\}$
and $\iota_{2}=\tilde{\varepsilon}_{r}+\ssum_{j=1}^{m-1}\varepsilon_{r,j}(\tilde{c}_{s,j})/\rho_{j+1}$.

By (\ref{eq:par_alpha}), (\ref{eq:e_r_bound}) in Theorem \ref{the:auxi_genera}
and $\textbf{DC}_{\xi}$, the bound of $\dot{V}_{m}$ 
becomes $\dot{V}_{m}\leq-\iota_{1}\alpha_{\xi}(\mu)V_{m}+\alpha_{\xi}(\mu)(\sigma\theta^{2}/2+\ssum_{j=1}^{m-1}\Xi_{j}/\rho_{j+1})
+\alpha_{\xi}(\mu)\iota_{2}\left(\gamma_{r}(\|e_{r}(t_{0})\|)\alpha(\mu(t_{0}))\right)^{2}\left(\alpha_{\xi}(\mu(t_{0}))\right)^{2L_{2}}$.
Let us choose any $\sigma$ satisfies (\ref{eq:sigma})
and $\bar{c}_{j}\geq\frac{1}{2}\sigma$ for $j=1,\cdots,m$, $\bar{\upsilon}_{j+1}\geq\frac{1}{2}\sigma$
for $j=1,\cdots,m-1$. As a result, the parameter $\tilde{c}_{m}$
is fixed and it is always possible to find an $h$ such that $\sigma\theta^{2}/(2\iota_{1})\leq h^{2}/8$. 
 Since $\Xi_{j}$ is independent of $\upsilon_{j+1}$ for $j=1,\cdots,m-1$,
it is independent of $\rho_{j+1}$. For any given $h$, it is always
possible to choose sufficiently large $\rho_{2},\cdots,\rho_{m}$
such that $\ssum_{j=1}^{m-1}\Xi_{j}/(\rho_{j+1}\iota_{1})\leq h^{2}/8$.
Then, $\tilde{\upsilon}_{m}$ is fixed. For any bounded initial condition
$e_{r}(t_{0})$ and $\mu(t_{0})$, it is always possible to find an
$h$ such that 
$
(\iota_{2}/\iota_{1})\left(\gamma_{r}(\|e_{r}(t_{0})\|)\right)^{2}\left(\alpha(\mu(t_{0}))\right)^{2}\left(\alpha_{\xi}(\mu(t_{0}))\right)^{2L_{2}}\leq h^{2}/4
$.
By  the upper inequality, one
has $h^{2}/2\geq\tilde{\epsilon}_{d}+\tilde{\epsilon}_{r}$ where
$\tilde{\epsilon}_{d}$ and $\tilde{\epsilon}_{r}$ are some finite constants.

When $V_{m}(\tilde{e}_{s}(t_{0}))\leq h^{2}/2,$ invoking comparison
lemma for $\dot V_{m}(\tilde{e}_{s})$ yields
$ V_{m}(\tilde{e}_{s}(t)) \leq\kappa^{-\iota_{1}}(\alpha_{\xi}(\mu))V_{m}(\tilde{e}_{s}(t_{0}))+(\tilde{\epsilon}_{d}+\tilde{\epsilon}_{r})(1-\kappa^{-\iota_{1}}(\alpha_{\xi}(\mu))) \leq h^{2}/2+\kappa^{-\iota_{1}}(\alpha_{\xi}(\mu))(V_{m}(\tilde{e}_{s}(t_{0}))-h^{2}/2)\leq h^{2}/2$,
which implies $\Omega(h)$ is an invariant set. It shows $\tilde{e}_{s}$
is bounded which implies that $\lim_{t\rightarrow t_{0}+T}e_{s}=0$
by (\ref{eq:norm_e_s-2}) in Lemma \ref{thm:6-2}. Similar to the
proof of Theorem \ref{thm:2}, we can prove that the DPTCO
problem is solved. \eproof
 
\begin{IEEEbiography}[{\includegraphics[width=1in, height=1.25in, clip, keepaspectratio]{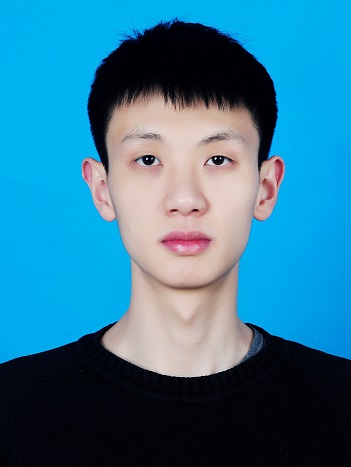}}]{Gewei Zuo}
	received  the M. E. Degree in control theory and engineering from Chongqing University, Chongqing, China, in 2022. He is currently pursuing the Ph.D. degree in control science and engineering with the school of Artificial Intelligence and Automation with Huazhong University of Science and Technology, Wuhan, Hubei, China.
	His research interests include Nonlinear System Control Theory, Distributed Cooperative Control and Distributed Convex Optimization.
\end{IEEEbiography}
\begin{IEEEbiography}[{\includegraphics[width=1in, height=1.25in, clip, keepaspectratio]{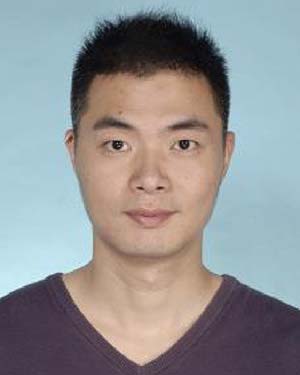}}]{Lijun Zhu}
	 received the Ph.D. degree in Electrical Engineering from University of Newcastle, Callaghan, Australia, in 2013. He is now a Professor in the School of Artificial Intelligence and Automation, Huazhong University of Science and Technology, Wuhan, China. Prior to this, he was a post-doctoral Fellow at the University of Hong Kong and the University of New- castle. His research interests include power systems, multi-agent systems and nonlinear systems analysis and control.
	\end{IEEEbiography}
\begin{IEEEbiography}[{\includegraphics[width=1in, height=1.25in, clip, keepaspectratio]{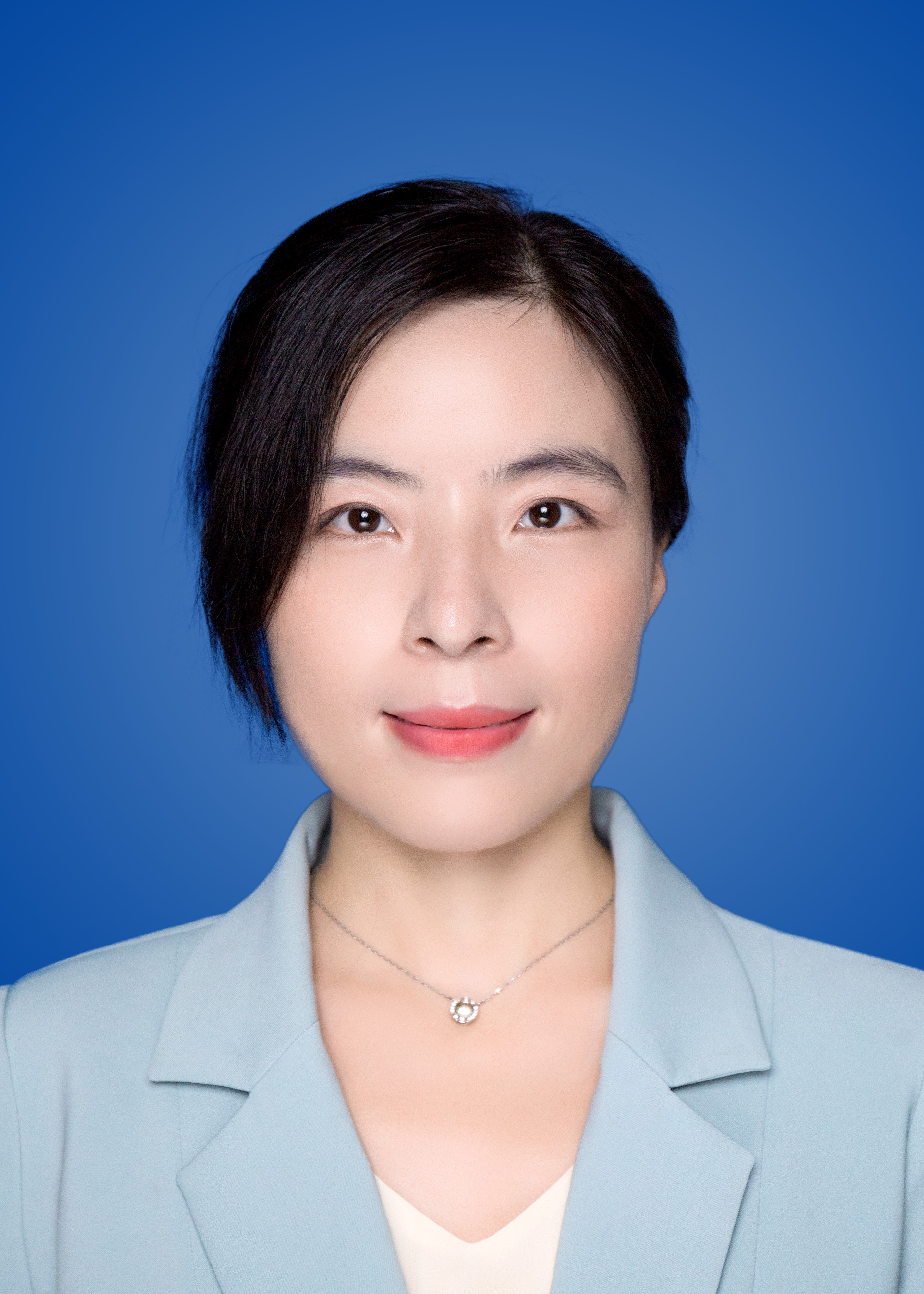}}]{Yujuan Wang}
	 received the Ph.D. degree in the School of Automation, Chongqing University, Chongqing, China, in 2016. She is now a Professor in the School of Automation, Chongqing University, Chongqing, China. Prior to this, she was a post-doctoral Fellow at the University of Hong Kong and a Joint Ph.D. Student at University of Texas at Arlington. Her research interests include nonlinear system control, distributed control, cooperative
adaptive control, fault-tolerant control.
	\end{IEEEbiography}
\begin{IEEEbiography}[{\includegraphics[width=1in, height=1.25in, clip, keepaspectratio]{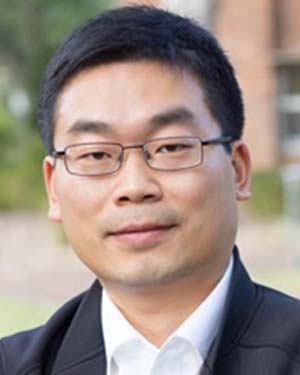}}]{Zhiyong Chen} (Senior Member, IEEE) received the B.E. degree in automation from the University
of Science and Technology of China, Hefei, China, in 2000 and the M.Phil. and Ph.D. degrees in mechanical and automation engineering from the Chinese University of Hong Kong, Hong Kong, in 2002 and 2005, respectively.
He was a Research Associate with the University of Virginia, Charlottesville, VA, USA, from
2005 to 2006. In 2006, he joined the University of Newcastle, Callaghan, NSW, Australia, where
he is currently a Professor. He was also a Changjiang Chair Professor with Central South University, Changsha, China. His research interests include nonlinear systems and control, biological systems, and reinforcement
learning.

Dr. Chen is/was an Associate Editor for \emph {Automatica, IEEE Transactions
on Automatic Control, IEEE Transactions on Neural
Networks and Learning Systems, and IEEE Transactions on Cybernetics}.
	\end{IEEEbiography}
\begin{IEEEbiography}
	[{\includegraphics[width=1in, height=1.25in, clip, keepaspectratio]{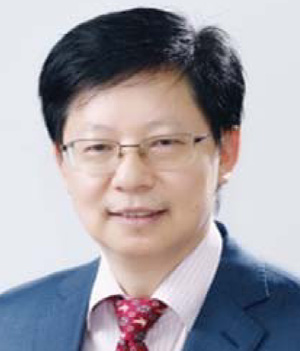}}]{Yongduan Song} (Fellow, IEEE) received the Ph.D. degree in electrical
	and computer engineering from Tennessee Technological
	University, Cookeville, TN, USA, in 1992. He held a
	tenured Full Professor with North Carolina A\&T State
	University, Greensboro, NC, USA, from 1993 to 2008
	and a Langley Distinguished Professor with the National
	Institute of Aerospace, Hampton, VA, USA, from 2005
	to 2008. He is currently the Dean of the Institute of
	Artificial Intelligence, Chongqing University, Chongqing,
	400044, China. He was one of the six Langley Distinguished
	Professors with the National Institute of
	Aerospace (NIA), Hampton, VA, USA, and the Founding Director of Cooperative
	Systems with NIA. His current research interests include intelligent systems,
	guidance navigation and control, bio-inspired adaptive and cooperative systems.
	Prof. Song was a recipient of several competitive research awards from the
	National Science Foundation, the National Aeronautics and Space Administration,
	the U.S. Air Force Office, the U.S. Army Research Office, and the U.S. Naval
	Research Office. He is an IEEE Fellow and has served/been serving as an
	Associate Editor for several prestigious international journals, including the \emph{IEEE
	Transactions on Automatic Control, IEEE Transactions on Neural Networks and
	Learning Systems, IEEE Transactions on Systems, Man, and Cybernetics: Systems,
	IEEE Transactions on Intelligent Transportation Systems}, etc. He is currently
	Editor-in-Chief for the \emph{IEEE Transactions on Neural Networks and Learning
	Systems}.
	\end{IEEEbiography}

\end{document}